\documentclass[12 pt]{amsart}
\usepackage{amscd, amssymb,url, tikz}
\RequirePackage{tikz}

\newtheorem{thm}{Theorem}

\newtheorem{lem}[thm]{Lemma}
\newtheorem{cor}[thm]{Corollary}
\theoremstyle{remark} % 'style changed again'
\newtheorem{rem}[thm]{Remark}% numbered with thm
\newtheorem{defn}[thm]{Definition}
\newtheorem{exa}[thm]{Example}% numbered with thm

\newcommand{\CC}{\mathbb C}

\newcommand{\NN}{\mathbb N}
\newcommand{\PP}{\mathbb P}
\newcommand{\QQ}{\mathbb Q}
\newcommand{\RR}{\mathbb R}

\newcommand{\ZZ}{\mathbb Z}

\newcommand{\Ocal}{\mathcal O}

\newcommand{\Hcal}{\mathcal H}

\newcommand{\rest}[1]{{\textstyle|}_{#1}}

\newcommand{\Aut}{\operatorname{Aut}}

\newcommand{\Diag}{\operatorname{Diag}}

\newcommand{\al}{\alpha}

\newcommand{\ga}{\gamma}

\newcommand{\la}{\lambda}

\newcommand{\cxi}{\mathtt{i}}

\newcommand{\trans}{{\mathsf T}}

\newcommand{\textsum}{{\textstyle{\sum}}}

\providecommand{\abs}[1]{\lvert#1\rvert}

\def\pmmu{{\pmb \mu}}

\def\RW{\operatorname{FJRW}}
\def\CR{\operatorname{CR}}

\def\wt{\widetilde}

\begin{document}
\title{LG/CY correspondence:\\ the state space isomorphism}
\author{Alessandro Chiodo}
\author{Yongbin Ruan}
\thanks{Partially supported by the National Science Foundation and
the Yangtze Center of Mathematics at Sichuan University}
%\date{}
\begin{abstract}
{We prove the  {\em classical mirror symmetry conjecture} for the mirror pairs
constructed by Berglund, H\"ubsch, and Krawitz. Our main tool is a
{\em cohomological LG/CY correspondence} which provides a degree-preserving isomorphism between the cohomology
of finite quotients of Calabi--Yau hypersurfaces inside a weighted
projective space and the Fan--Jarvis--Ruan--Witten state space of
the associated Landau--Ginzburg singularity theory. }
\end{abstract}

\maketitle

\pagestyle{plain}

\section{Introduction}
Mirror symmetry has been one of the most inspirational problems
arising from physics in the last twenty years. In the most common
formulation, which we call \emph{classical mirror symmetry}, it is
a duality statement pairing two Calabi--Yau three-folds
$X^3$ and $Y^3$ by interchanging $h^{1,1}$ and $h^{2,1}$.
When the mirror symmetry was first proposed twenty years ago, only a few examples of Calabi--Yau
three-folds were known. A major effort was launched to construct more examples.
Soon, physicists constructed millions of examples which are
(orbifolded) hypersurfaces and complete intersections lying inside
weighted projective spaces or toric varieties. Since every three-dimensional Calabi--Yau
orbifold admits a crepant resolution, we obtain millions of examples of smooth
Calabi--Yau three-folds.

Among these millions of examples, an elementary and yet elegant mirror symmetry
construction was proposed by
the physicists Berglund and H\"ubsch \cite{BH}, which will be the
focus of our interest. In \cite{BH} a hypersurface $X_W$ in a
weighted projective space $\PP(\pmb w)=\PP(w_1, \dots, w_N)$ is
considered: $X_W$ is defined by a quasihomogeneous polynomial $W$.
Berglund and H\"ubsch describe a simple definition of the mirror
of $X_W$.

The construction only involves cases when $W$ is ``invertible'';
\emph{i.e.} $W$ is the sum of $N$ monomials, as many as the variables.
In this case, one can transpose the
exponents matrix and obtain another quasihomogeneous
polynomial $W^\trans$ defining a hypersurface lying
in another weighted projective space. The varieties
$X_W$ and $X_{W^T}$ are not mirror pairs in general and a certain
orbifolding construction must be involved.
Berglund and H\"ubsch proposed a certain physical property for  a correspondence between automorphism groups
$G\subset \Aut(\{W=0\})$ and $G^{\trans}\subset\Aut(\{W^\trans=0\})$; in \cite{BH},
the Calabi--Yau $X_W/G$
is expected to be the
mirror image of the Calabi--Yau $X_{W^\trans}/G^\trans$.
More precisely, the classical mirror symmetry conjecture should hold for these pairs:
if we stick to Calabi--Yau three-folds,
$h^{1,1}$ and $h^{2,1}$ should be interchanged.
This group duality is precisely stated only in some cases, but
already opens the way to several interesting tests:
Kreuzer and Skarke checked  thousands three-folds for which
they computed the so-called ``Landau--Ginzburg phase'' \cite{KStest}. Indeed, these invariants
exhibit the classical mirror symmetry correspondence.

Unfortunately, this approach was mysteriously  abandoned  to favor a more geometric approach due to Batyrev and Borisov.
Batyrev and Borisov considered the complete intersection of a  Gorenstein toric variety. In this context,
the mirror symmetry was interpreted as polar duality. A major theorem of the day was a solution of
the classical mirror symmetry conjecture in this context. We should mention that Batyrev--Borisov imposed an important
condition called {\em Gorenstein} in all their constructions. Indeed, Gorenstein conditions
are also crucial on our recent
investigation of Gromov-Witten theory \cite{CIR}.
It is interesting to consider it in
the context of weighted projective spaces. The ambient weighted projective space $\PP(\pmb w)$ is
Gorenstein if and only if $\sum_j w_j$ is a multiple of every
weight $w_j$; hence, with a Gorenstein ambient space we can reduce
to the Calabi--Yau hypersurface defined by the Fermat polynomial
of degree $d=\sum_jw_j$; \emph{i.e.} $W(x_1,\dots,x_N)=\sum_j
x^{d/w_j}_j$. It was known that Fermat Calabi--Yau hypersurfaces only represents a small subclass of
all Calabi--Yau hypersurfaces. It was a big surprise to us that  a vast range of cases involved in the
Berglund--H\"ubsch construction are not covered by Batyrev and
Borisov (see Remark \ref{rem:toric})!

During the last two years, interest in this problem was revived
by the introduction of
a Gromov--Witten-type theory for singularities by Fan, Jarvis, and the second author.
This fits within the framework of
the Landau--Ginzburg  model and is based
on a proposal of Witten (FJRW theory).
Recently, Krawitz \cite{Kr} found
a general construction for the dual group $G^\trans$. Working on
much more general grounds where $X_W$ is not necessarily
Calabi--Yau, Krawitz proved an ``LG-to-LG'' mirror symmetry
theorem for all invertible polynomials $W$ and all
admissible groups $G$.
We should emphasize that the Berglund--H\"ubsch--Krawitz computations
are purely in the Landau--Ginzburg setting. Whether $X_W/G$ and
$X_{W^\trans}/G^\trans$ are a mirror pair of Calabi--Yau orbifolds
is an open question.
We shall give a firm answer in this article. To state our theorem, let us set up some notation.

\smallskip

\noindent \textbf{The mirror symmetry setup.} A hypersurface inside a
weighted projective space is defined by a
quasihomogeneous polynomial $W$
in the variables $x_1,\dots,x_N$ of charges
$q_1, \dots,q_N\in \QQ_{>0}$
such that
\begin{equation}\label{eq:qhom}
W(\lambda^{q_1}x_1, \dots, \lambda^{q_N}x_N)=\la W(x_1, \dots, x_N).\end{equation}
Write $q_1=w_1/d, \dots, q_N=w_N/d$ with common denominator
so that we have  $\gcd(w_1, \dots, w_N, d)=1$.
Then, $X_W=\{W=0\}\subset \PP(w_1, \dots, w_N)$
defines a degree-$d$ hypersurface. We always assume that
$W$ has a unique singularity at zero; in other words, $X_W$
is a smooth Deligne--Mumford stack (an orbifold).
Furthermore, $X_W$ is a Calabi--Yau orbifold
if and only if $\sum_j q_j=1$; we refer to this condition as  the \emph{CY condition}
(see also Section \ref{sect:CYside}).
For three-dimensional Calabi--Yau orbifolds, the
crepant resolution always exists and the Hodge numbers are equal
to the Hodge numbers of the underlying Chen--Ruan orbifold
cohomology.
A wider range of Calabi--Yau orbifolds
arises from quotients of $X_W$.
Consider the group $\Aut(W)$ of diagonal symmetries
rescaling the coordinates and preserving $W$:
$(\al_1,\dots,\al_N)\in \CC^\times$ such that
$W(\al_1 x_1,\dots,\al_Nx_N)$ equals $W(x_1,\dots,x_N)$
for all $(x_1,\dots,x_N)\in \CC^N.$
Clearly $J_W=(\exp(2\pi \cxi q_1),
 \dots, \exp(2\pi \cxi q_N))$ is contained in $\Aut(W)$ and
 the action of $J_W$ on $X_W$ is trivial (see Section \ref{sect:CYside} for
 a discussion of group actions on these stacks).
 For any subgroup $G$ of diagonal symmetries containing $J_W$, let
 us consider the group $\wt{G}=G/\langle J_W\rangle$ acting faithfully on $X_W$.
 The quotient is Calabi--Yau as long as $G$ is contained in $SL_N(\CC)$.
 Let $G\subset \Aut(W)$ be such that
   $\langle J_W\rangle \subseteq G\subseteq SL(N,\CC)$. Then,
   there is a very natural construction associating
   $W^\trans$ and $G^\trans$ to $W$ and $G$ and preserving the
   following properties
   (see Section \S\ref{sect:BH}, \eqref{eq:daulpoly} and \eqref{eq:dualgroup}, for a concise
   presentation of the construction of
   $W^\trans$ and $G^\trans$).

   First, the polynomial $W^\trans$ --- precisely as the polynomial $W$ --- the polynomial
   $W^\trans\colon \CC^N\to \CC$ has a unique singularity at $\pmb 0$ and the sum of its charges
$q_1^\trans,\dots, q_N^\trans$ equals $1$
(\emph{i.e.} $X_{W^\trans}$ is Calabi--Yau).
Second, the group $G^\trans$ --- in perfect analogy with $\langle J_W\rangle \subseteq G\subseteq SL(N,\CC)$ ---
satisfies $\langle J_{W^\trans}\rangle
\subseteq G^\trans\subseteq SL(N,\CC)$.

\medskip
Our mirror symmetry theorem is

\smallskip

\noindent\textbf{Theorem \ref{cor:MS}}.
\emph{
The Calabi--Yau $[X_W/\wt G]$ and
the Calabi--Yau
$[X_{W^\trans}/\wt G^\trans]$ form a mirror pair; \emph{i.e.} we have
$$H_{\CR}^{p,q}([X_W/\wt G];\CC)\cong H_{\CR}^{N-2-p,q}([X_{W^\trans}/\wt G^\trans];\CC),$$
where $H_{\CR}(\ \ ;\CC)$ stands for Chen--Ruan orbifold cohomology.}

\smallskip

   The above theorem is  precisely performing a ``90 degrees rotation of
   the Hodge diamond'' as predicted by
   the classical mirror symmetry conjecture in these cases.

 \begin{rem}Let us point out that one can find two different polynomials $W_1, W_2$
   in the same family of degree-$d$ quasihomogeneous polynomials in the variables
   $x_1,\dots,x_N$ with charges $q_1,\dots,q_N$. Now, whereas $X_{W_1}$ may be regarded as a
   deformation of $X_{W_2}$,
   there is no apparent reason to claim that
   $W_1^\trans$ is related to $W_2^\trans$. Indeed the above
   statement implies that the cohomologies of the hypersurfaces
   defined by $W_1^\trans$ and $W_2^\trans$ are strictly related (in many cases, \emph{e.g.}
   when $SL_W=\langle J_W\rangle$, they are isomorphic).
   This provides many examples of ``multiple
   mirrors'' which are not birational to each other---a
   rather interesting phenomenon which is certainly
   worth further investigation.
   \end{rem}

   \noindent {\bf The LG/CY correspondence.} The above mirror symmetry theorem is an outcome of our program to study
   so called {\em Landau--Ginzburg (LG)/Calabi--Yau (CY) correspondence}.
    In the early days of mirror symmetry, physicists
    noticed that regarding $W$ as a function on $\CC^N$ leads to the
    Landau--Ginzburg (LG) singularity model. (In this correspondence,
    we place ourselves within a more general framework: we do not need to
    require that the number of variables equals the number of monomials.)
    The argument has been made on physical
    grounds  that there should be a
    LG/CY correspondence connecting
    Calabi--Yau  geometry to the
    LG singularity model \cite{VW} \cite{Wi3}. In this context,
    CY manifolds are considered from the point of view of
    Gromov--Witten theory; this correspondence would therefore inevitably
    yield new predictions on Gromov--Witten invariants and is likely
    to greatly simplify their calculation (it is generally believed that
    the LG singularity model is relatively easy to compute).
    In a different context, the LG/CY correspondence led to identifying
    matrix factorization as the LG counterpart of the derived category
    of complexes of coherent
    sheaves \cite{H}, \cite{Ko1}.

    In \cite{FJR1, FJR2, FJR3}, a candidate quantum theory of singularities
    has been constructed by Fan, Jarvis, and Ruan.
        Using the Fan--Jarvis--Ruan--Witten theory as a candidate theory on the  LG side,
    the authors have launched a program to solve
    LG/CY-correspondence for Calabi--Yau hypersurfaces inside
    weighted projective spaces. In \cite{ChRu},
    the equivalence between FJRW theory and GW theory
    has been established in genus zero in the case of
    the famous quintic three-fold. The starting point of this equivalence is
    an isomorphism between the two cornerstones the two theories are built upon:
    the $\RW$ state space of the singularity and the cohomology of the hypersurface.
    This can be done explicitly in several examples, but
    it is rather intricate to prove it in full generality (see Section \ref{sect:warmup}
    for a case-by-case approach through elliptic curves, K3 surfaces and
    Calabi--Yau three-folds).

    We will
    accomplish the isomorphism in full generality by building a common combinatorial model for both
    theories. Our model generalizes the
    combinatorial model of Boissi\`ere, Mann and Perroni \cite{BMP} for weighted projective spaces.
    The main result
is the following {\em cohomological LG/CY correspondence} where
$H_{\CR}^{p,q}([X_W/\wt G];\CC)$ denotes the Chen--Ruan orbifold cohomology while $H^{p,q}_{\RW}(W,G;\CC)$
denotes the state space of Fan-Jarvis-Ruan-Witten theory (see Section \ref{sect:sides} for the detailed definition).

\medskip

\noindent{\textbf{Theorem \ref{thm:main}.}
\emph{Let $W$
be a nondegenerate
quasihomogeneous polynomial of degree $d$ in the variables $x_1,\dots,x_N$
whose charges add up to $1$ (CY condition).
Then, for any group $G$ of diagonal symmetries containing $J_W$
we have a bidegree preserving isomorphism of vector spaces
$$H_{\CR}^{p,q}([X_W/\wt G];\CC)\cong H^{p,q}_{\RW}(W,G;\CC).$$}

The mirror symmetry theorem is a direct consequence of our cohomological LG/CY correspondence and Krawitz
LG-to-LG mirror symmetry theorem.

We point out, however, a most surprising aspect of our main theorem:
not only does it hold for noninvertible polynomials, it also
holds for $G\not \subseteq SL_W$ (\emph{e.g.} $G$ equal to
the group $\Aut(W)$ itself).
This goes beyond the
LG/CY-correspondence stated in physics and yields several surprising
consequences.

   \subsection{Structure of the paper}
   This article is organized as follows.
   In Section \ref{sect:BH} we state precisely the mirror symmetry
   construction.
   In Section \ref{sect:sides},
   we introduce the
   state spaces of both Gromov--Witten theory
   (CY side) and Fan--Jarvis--Ruan--Witten theory (LG side) and
   we state the cohomological
   Landau--Ginzburg(LG)/Calabi--Yau(CY) correspondence between them.
   In Section \ref{sect:warmup} we
   present several examples illustrating the correspondence, this prepares the ground
   to the combinatorics involved in the general proof.
   In Section \ref{sect:main}, we prove the two theorems stated above.
   In Section \ref{sect:exa} we review the examples introduced in Section \ref{sect:warmup}
   in the light of the combinatorial tools introduced in Section \ref{sect:main}.

\section{The classical mirror symmetry construction}\label{sect:BH}
Berglund and H\"ubsch \cite{BH}
consider
polynomials in $N$ variables having $N$ monomials
\begin{equation}\label{eq:invertible}
W(x_1,\dots,x_N)=\sum_{i=1}^N\prod_{j=1}^N x_j^{m_{i,j}}.
\end{equation}
Note that each of the $N$ monomials
has coefficient one; indeed,
since the number of variables equals the number of monomials,
even when we start from a polynomial of the form
$\sum_{i=1}^N l_i \prod_{j=1}^N x_j^{m_{i,j}}$, it is possible to reduce to
the above form by conveniently
rescaling the $N$ variables.
In this way assigning a polynomial $W$ as above amounts
to specifying a square matrix
$$M=(m_{i,j})_{1\le i,j\le N}.$$

The polynomials studied in \cite{BH} are called ``invertible'', because the matrix
$M$ is an invertible $N\times N$ matrix. In fact, polynomials of this type
may be regarded as quasihomogeneous polynomials in the variables
$x_1,\dots,x_N$ of charges $q_1,\dots,q_N$ as in \eqref{eq:qhom}: to this effect,
simply set
\begin{equation}\label{eq:chargesfromM}q_i=\textsum_i m^{i,j},
\end{equation}
the sum of the entries on the $i$th line of $M^{-1}=(m^{i,j})_{1\le i,j\le N}$.

Each column $(m^{1,j},\dots,m^{N,j})^\trans$
of the matrix $M^{-1}$ can be used to
define the diagonal matrix
$\rho_j$ whose diagonal entries are $\exp(2\pi \cxi m^{1,j}), \dots,$
and $\exp(2\pi \cxi m^{N,j})$. In fact these matrices satisfy the following property
$\rho_j^*W=W$; \emph{i.e.} $W$ is invariant with respect to
$\rho_j$.
Furthermore the group
$\Aut(W)$ of diagonal matrices $\al$ such that $\al^*W=W$ is generated
by the elements $\rho_1,\dots,\rho_N$:
$$\Aut(W):=\{\al={\rm Diag}(\al_1,\dots,\al_N)\mid \al^*W=W\}=
\langle \rho_1,\dots,\rho_N\rangle.$$
For instance, the above mentioned matrix $J_W$ whose diagonal entries are $\exp(2\pi\cxi q_1),\dots,$ and
$\exp(2\pi\cxi q_N)$
lies in $\Aut(W)$ and is indeed the product $\rho_1\cdots\rho_N$.
Write $$SL_W=\Aut(W)\cap {\rm SL}_\CC(N),$$
the matrices with determinant $1$;
in Berglund and H\"ubsch's construction we consider groups $G$
containing $J_W$ and contained in $SL_W$. We write $\wt G$ for the
quotient $G/\langle J_W\rangle$.

Indeed, we start form a polynomial $W\colon \CC^N\to \CC$
of nondegenerate type: \emph{i.e.} having
a single critical point, the origin $\pmb 0\in \CC^N$.
Let $W$ be a nondegenerate invertible potential
of charges $q_1,\dots, q_N$ satisfying the \emph{Calabi--Yau condition}
\begin{equation}\label{eq:CYsect2condition}\sum_j q_j=1.
\end{equation}
The geometrical meaning of this condition is the following:
$X_W=\{W=0\}$ is Calabi--Yau, or --- more precisely --- $\{W=0\}$ is a degree-$d$
{Calabi--Yau} hypersurface in the weighted projective space
$\PP(d\pmb q)$, where $d$ is the least integer for which  $d\pmb q\in \ZZ^N$.
Let $G\subset \Aut(W)$ be a group of diagonal symmetries
satisfying $\langle J_W\rangle \subseteq G\subseteq SL_W$ (the fact that $J_W$ is contained
in $SL_W$ follows from the Calabi--Yau condition).

In this context there is a natural way to associate to $W$ a
polynomial $W^\trans$ and to $G$ a subgroup $G^\trans$ of
the group of diagonal symmetries of the polynomial $W^\trans$.
The polynomial $W^\trans$ is defined by
transposing the matrix $(m_{i,j})$:
\begin{equation}\label{eq:daulpoly}
W^\trans(x_1,\dots,x_N)=\sum_{i=1}^N\prod_{j=1}^N x_j^{m_{j,i}}.
\end{equation}
The group $G^\trans$ is defined by
\begin{equation}\label{eq:dualgroup}
G^\trans=\textstyle{\left \{\prod_{j=1}^N
(\rho^\trans_i)^{a_i} \mid \text{ if }\prod_{j=1}^N
x_i^{a_i} \text{ is } G\text{-invariant}\right\}},
\end{equation}
where $\rho^{\trans}_i$ is the diagonal symmetry corresponding to
the $i$th column of $(M_{W^\trans})^{-1}$ (note that, by
construction, $M_{W^\trans}$ equals $(M_{W})^\trans$).

Then we
have the following properties:
\begin{itemize}
\item[--]
$W^\trans$ is nondegenerate and the sum of its charges
$q_1^\trans,\dots, q_N^\trans$ equals $1$
(\emph{i.e.} $X_{W^\trans}$ is Calabi--Yau).
\item[--]
The group $G^\trans$ satisfies $\langle J_{W^\trans}\rangle
\subseteq G^\trans\subseteq SL_{W^\trans}$.
\item[--]
The quotients $[X_W/\wt G]$ and $[X_{W^\trans}/\wt G^\trans]$
form a mirror pair in the following
sense.
\end{itemize}

Below,
$M$, $W$, and $G$ satisfy these conditions:
$M=(m_{i,j})$ is an invertible $N\times N$ matrix satisfying
$\sum_{i,j} m^{i,j}=1$ (CY condition, see \eqref{eq:chargesfromM} and \eqref{eq:CYsect2condition}),
the polynomial $W(x_1,\dots, x_N)=\sum_i \prod_j x_j^{m_{i,j}}$ has
a single isolated critical point at $\pmb 0\in \CC^N$,
$G$ is a group containing $J_W$ and contained in $SL_W$.

\begin{thm}\label{cor:MS}
Then, the Calabi--Yau $[X_W/\wt G]$ and
the Calabi--Yau
$[X_{W^\trans}/\wt G^\trans]$ form a mirror pair; \emph{i.e.} we have
$$H_{\CR}^{p,q}([X_W/\wt G];\CC)\cong H_{\CR}^{N-2-p,q}([X_{W^\trans}/\wt G^\trans];\CC),$$
where $H_{\CR}(\ \ ;\CC)$ stands for Chen--Ruan orbifold cohomology.
\end{thm}

We prove this theorem in Section \ref{sect:main}.

\begin{rem}
Let us mention that the fact that $W^\trans$ is nondegenerate follows from
Kreuzer and Skarke \cite{KS} classification
of invertible nondegenerate
potentials. An invertible potential $W$ is nondegenerate if and only if it can be written as a
sum of (decoupled) invertible potentials of one of the following
three types, which we will refer to as \emph{atomic types}:
\begin{itemize}
\item[] $W_\text{Fermat} = x^a.$
\item[] $W_\text{loop}= x_1^{a_1}x_2+x_2^{a_2}x_3+\dots +x_{N-1}^{a_{N-1}}x_N+x_N^{a_N}x_1.$
\item[] $W_\text{chain}= x_1^{a_1}x_2+x_2^{a_2}x_3+\dots +x_{N-1}^{a_{N-1}}x_N+x_N^{a_N}.$
\end{itemize}
If $W$ is a Fermat type polynomial, the
ambient weighted projective stack is Gorenstein. However, if $W$
is of a loop or chain type, the ambient weighted projective stack
is not Gorenstein in general.
\end{rem}

\begin{cor}\label{cor:coarsecohom}
Assume that the quotient schemes $X_W/\wt G$ and
$X_{W^\trans}/\wt{G}^\trans$ both admit crepant resolutions $Z$
and $Z^\trans$. Then the above statement yields a statement in
ordinary cohomology:
$$h^{p,q}(Z;\CC)=h^{N-2-p,q}(Z^\trans;\CC).$$
\end{cor}

We prove this corollary in Section \ref{sect:main}.

\begin{rem}\label{rem:toric} In the case where $w_j$ divides $d$,
Theorem \ref{cor:MS} can be deduced from Batyrev's
construction of mirror pairs into toric geometry. The general
case does not fit in this framework because the ambient space (unlike
the space $X_W$) is not
Gorenstein in general. The following example illustrates this well.
\end{rem}

\begin{exa}
In order to illustrate the above statement we provide an example straight-away and
we refer to Section \ref{sect:exa} for more.
Consider the following quintic
hypersurface in $\PP^4$:
$$\{x_1^4x_2+x_2^4x_3+x_3^4x_4+x_4^4x_5+x_5^5=0\}_{\PP^4}.$$
It is a chain-type Calabi--Yau variety $X$ whose Hodge diamond is clearly equal to
that of the Fermat quintic and is well known: $h^{1,1}=1$, $h^{0,3}=1$,
$h_{1,2}=101$
\begin{equation}\begin{matrix}\label{eq:quinticdiamond}
& & &   &1&   &   \\
& & &0  & &0  &   \\
& &0&   &1&   &0  \\
&1& &101& &101& &1\\
& &0&   &1&   &0  \\
& & &0  & &0  &   \\
\quad & & &   &1&   &&&.
\end{matrix}\end{equation}
The mirror Calabi--Yau
is given by
the vanishing of the polynomial
$$W^\trans(x_1,x_2,x_3,x_4,x_5)=x_1^4+x_1x_2^4+x_2x_3^4+x_3x_4^4+x_4x_5^5=0,$$
which may be regarded as defining a degree-$256$
hypersurface $X^\trans$ inside
$\PP(64,48,52,51,41)$. This degree-$256$ hypersurface is Calabi--Yau
(\emph{i.e.} $256$ is indeed the sum of the weights). But the
ambient weighted projective stack is no longer Gorenstein. Note
that the group $SL_{W^\trans}$ coincides with $\langle
J_W^\trans\rangle$; therefore Corollary \ref{cor:MS} reads
$$h^{p,q}(X;\CC)=h^{3-p,q}(X^\trans;\CC).$$
Indeed, the Hodge diamond satisfies
$h^{1,1}=101$, $h^{0,3}=1$,
$h_{1,2}=1$.
$$\begin{matrix}
 & &   &1&   &   \\
 & &0  & &0  &   \\
 &0&   &101&   &0  \\
\ \ \ \ \ \ 1& &  1&   &1  & &1\\
 &0&   &101&   &0  \\
 & &0  & &0  &   \\
\hspace{1cm} & &   &1&   &&&.
\end{matrix}$$

%In cases like the one illustrated above one may envisage choosing
%a different embedding and embedding also the three-fold $X^\trans$
%into Gorenstein toric varieties. This has been done in special
%cases (See \cite{CDGP} and references therein for some
%calculations in special cases). However, the authors are not aware
%of any general procedure.
\end{exa}

\section{The cohomological LG/CY correspondence}\label{sect:sides}
The geometrical Landau--Ginzburg/Calabi--Yau correspondence is a
correspondence between two geometrical settings defined starting
from the polynomial $W$ and the group $G$. With respect to the
previous section we work in a more general setup.

\subsection{The polynomial and its diagonal symmetries.}
We consider polynomials
\begin{equation}\label{eq:Wexplicit}
W(x_1,\dots,x_N)=l_1 \prod_{j=1}^N x_j^{m_{1,j}}+\dots
+l_s \prod_{j=1}^N x_j^{m_{s,j}}.
\end{equation}
where $l_1, \dots, l_s$ are nonzero complex numbers
and $m_{i,j}$ (for $1\le i \le N$ and $1\le j\le s$)
are nonnegative integers. We will always suppose that the summands of the above decomposition
are distinct monomials; \emph{i.e.} monomials with distinct exponents.

We assume that $W$ is \emph{quasihomogeneous}; \emph{i.e.} there
exist positive integers $w_1,\dots, w_N$, and $d$ satisfying
\begin{equation}\label{eq:quasihomog}
W(\lambda^{w_1}x_1, \dots, \lambda^{w_N}x_N)=\lambda^{d}
W(x_1, \dots, x_N) \qquad \forall \la \in \CC,\end{equation}
or, equivalently,
\begin{equation}\label{eq:quasihomogdiff}
W=\textstyle{\sum_{j=1}^N \frac{w_j}d x_j {\partial_jW}}
\end{equation}
(we write $\partial_j$ for the
partial derivative with respect to the $j$th variable).
For $1\le j\le N$, we say that the \emph{charge} of the variable
$x_j$ is $q_j=w_j/d$. As soon as $w_1, \dots, w_N$ and $d$ are
coprime, we say that the \emph{degree} of $W$ is $d$ and that the \emph{weight} of the variable $x_j$ is $w_j$.
We assume
that the origin is the only \emph {critical point} of $W$;
\emph{i.e.} the only solution of
\begin{equation}\label{eq:critical}
{\partial_jW}(x_1,\dots,x_N)=0 \qquad \text{for
} j=1,\dots, N\end{equation} is $(x_1,\dots, x_N)=(0,\dots,0)$.
(By \eqref{eq:quasihomogdiff}, if $(x_1,\dots,x_N)$
satisfies \eqref{eq:critical}, then $W(x_1,\dots, x_N)$ is zero.)

\begin{defn}\label{defn:qhomsing}
We say that $W$ is a \emph{nondegenerate quasihomogeneous
polynomial} if it is a quasihomogeneous polynomial of degree $d$
in the variables $x_1, \dots, x_N$ of charges $w_1/d,\dots, w_N/d>0$
and the following conditions are satisfied:
\begin{enumerate}
\item
$W$ has a single critical point at the origin;
\item \label{chargesunique}
the charges are uniquely determined by $W$.
\end{enumerate}
\end{defn}

\begin{rem}\label{rem:rank}
The second condition above may be regarded as saying that the
$s\times N$ matrix $M_W=(m_{i,j})$
defined by $W(\pmb x)=\sum_{i=1}^s
l_i \prod_{j=1}^N x_j^{m_{i,j}}$ has rank $N$ (\emph{i.e.} has a
left inverse).
\end{rem}

\noindent \emph{CY condition.} The main result of this paper,
the cohomological Landau--Ginzburg/Calabi--Yau correspondence, holds
under the following condition:
\begin{equation}\label{eq:CYconditionSectLG/CY}
\sum_j q_j=1.\end{equation}

The definition of $\Aut(W)$ applies without changes
to the polynomial $W$ in
this context: $\Aut(W)$ is
the group of $(\al_1,\dots,\al_N)\in (\CC^\times)^N$
satisfying $W(\al_1x_1,\dots,\al_Nx_N)=
W(x_1,\dots,x_N)$. Again  $SL_W=SL(\CC,N)\cap \Aut(W)$
and $J_W:=(\exp({2\pi {\cxi}{w_1}/{d}}), \dots, \exp({2\pi
\cxi{w_N}/{d}}))$ is in $SL_W$
and generates a cyclic subgroup of order $d$ as
a consequence of the CY condition.

\subsection{The Calabi--Yau side.}\label{sect:CYside}
On the Calabi--Yau side the picture is that of a hypersurface inside the
weighted projective stack\footnote{From now on we will always stress the stack-theoretic
nature of the above quotient, because this point of view is crucial here.}
$$\PP(w_1,\dots,w_N)=[(\CC^N\setminus\{\pmb 0\})/\CC^\times],$$
where $\CC^\times$ acts as $\la(x_1,\dots,x_N)=(\la^{w_1}x_N,\dots,\la^{w_N}x_N)$ and
$w_1,\dots,w_N$ are the weights satisfying $q_j=w_j/d$.
By the nondegeneracy condition,
the equation $W=0$ defines a smooth hypersurface inside
$\CC^N\setminus \{\pmb 0\}$:
the normal vector
\begin{equation}\label{eq:normalvector}\vec{n}(\pmb x)=\left(\partial_j W(\pmb x)\right)_{j=1}^N\end{equation}
never vanishes on $\CC^N\setminus \{\pmb 0\}$.
By the quasihomogeneity condition the action of $\CC^\times$ fixes the variety $\{W=0\}$.
We write $X_W$ for the quotient stack
$$X_W:=[\{W=0\}_{\CC^N\setminus \{\pmb 0\}}/\CC^\times]\subset \PP(w_1,\dots,w_N).$$

\begin{rem} \label{rem:CY}
Note that the CY condition implies that $\omega_{X_W}$ is trivial,
$X_W$ has canonical singularities,
and $H^i(X_W,\Ocal_{X_W}) = (0)$ for $i = 1, \dots, n-1$ (see \cite[Lem.~1.12]{CG}).
In other words
$X_W$ is Calabi--Yau (see \cite[4.1.8]{Bat94}).
We point out that well-formedness conditions (see \cite{IanoF} and \cite[p.8]{CG})
are not needed here, see Remark \ref{rem:norbehaviour} in Section \ref{sect:main}.
\end{rem}

Consider a group $G$ contained in $\Aut(W)$ and containing $J_W$. The homomorphism
mapping $\la \in \CC^\times$ to
$(\la^{w_1},\dots,\la^{w_N})\in (\CC^\times)^N$, is injective because $\cap _j
\pmmu_{w_j}$ is trivial (the weights are coprime by definition).
It is natural to identify $\CC^\times$ with the image of the above
injection: we write $\bar \la$ for the image of $\la\in
\CC^\times$, \emph{i.e.}
\begin{equation}\label{eq:barla}
\bar \la =(\la^{w_j})_{j=1}^N.
\end{equation}
Notice that we have
\begin{equation}\label{eq:intgroups}
\CC^\times \cap G=\langle J_W\rangle\end{equation}
as a straightforward consequence
of the quasihomogeneity of $W$.
The group $\wt G=G/\langle J_W\rangle$ acts faithfully on the stack $X_W$.
In fact, following Romagny's treatment \cite{Ro} of actions on stacks
we may consider the $2$-stack $[X_W/\wt G]$ which is equivalent to
the quotient stack of $\{W=0\}_{\CC^N\setminus \{\pmb 0\}}$
by the action of the product
$$G\CC^\times=\{g(\la^{w_1},\dots,\la^{w_N})\mid g \in G\subset (\CC^\times)^N,\la \in \CC^\times\}
\subseteq(\CC^\times)^N$$
(this is a consequence of $G\CC^\times /\CC^\times=\wt G$ and of
 \cite[Rem.~2.4]{Ro}). In this way
we may exhibit $[X_W/\wt G]$ as a  quotient and indeed a
smooth stack of Deligne--Mumford type:
\begin{equation}\label{eq:quotstack}
[X_W/\wt G]=[\{W=0\}_{\CC^N\setminus \{\pmb 0\}}/G\CC^\times] \qquad
(\text{with $\wt G=G/\langle
J_W\rangle$}).
\end{equation}
Alternatively, one may take the above formula as a definition of
the quotient $[X_W/\wt G]$.

\begin{rem} If $G\subseteq SL_W$, the $G$-action preserves the canonical
form on $X_W$ and the quotient space $Y=X_W/\wt G$
is still Calabi--Yau (see Remark \ref{rem:norbehaviour} in
Section \ref{sect:main}).
This motivates the hypothesis
$G\subseteq SL_W$ in \cite{BH}; however, Theorem \ref{thm:main}
holds for the orbifold $[X_W/\wt G]$ even beyond $SL_W$.
This happens because the theorem is phrased in terms of Chen--Ruan
orbifold cohomology and applies to an orbifold which
--- in some sense --- is Calabi--Yau
(the CY condition insures that the canonical divisor $K$ of the stack $[X_W/G]$
has vanishing degree).
Example \ref{exa:groupquot} exhibits a situation where $G=\Aut(W)$ is
not contained in $SL_W$; there, the quotient space $X_W/\wt G$
is not Calabi--Yau (it is a projective line)
but the stack $[X_W/\wt G]$ has canonical divisor of degree $0$ and
in fact there exists a tensor power of the canonical line bundle which is trivial
(we have $\omega^{\otimes 4}\cong \Ocal$). This
is enough for Theorem \ref{thm:main} on LG/CY correspondence
to hold at a stack-theoretic level even if there is
no scheme-theoretic counterpart to this statement.
\end{rem}

The main invariant on the Calabi--Yau side is
the \emph{Chen--Ruan orbifold cohomology}.
For a smooth Deligne--Mumford quotient stack $\mathcal X=[U/G]$
it may be regarded essentially as follows. It is a
direct sum over the elements $g$ of the group $G$: the summands are ordinary
cohomology groups $H^{\bullet}(\cdot \ ;\CC)$ of the so-called sectors $\mathcal X_g=[\{u\in U\mid gu=u\}/G]$.
The sectors are algebraic stacks of Deligne--Mumford type; since the cohomology with complex coefficients
can be identified with the cohomology of the coarse space, the summands
can be expressed in
terms of coarse spaces. We now detail this description for the quotient stack
$$[X_W/\wt G]=[\{W=0\}_{\CC^N\setminus \{\pmb 0\}}/G\CC^\times].$$

For any $\ga\in (\CC^\times)^N$, and in particular for $\ga \in G\CC^\times$,
we can define
\begin{align}\label{eq:fixlocus}
\CC^N_{\gamma}&= \{\pmb x\in \CC^N \mid \gamma \pmb x=\pmb x\};\\
N_\gamma&=\dim_{\CC}(\CC^N_\gamma) \label{eq:fixdim};\\
W_{\gamma}&=W\rest{\CC^N_{\gamma}}\label{eq:restfix}.
\end{align}
For $\ga \in G\CC^\times$,
we set the notation $$\{W_\ga=0\}_\ga:=\{W_\ga=0\}_{\CC^N_\ga\setminus \{\pmb 0\}};$$
it is easy to show that $\{W_\ga=0\}$ defines a smooth hypersurface inside $\CC_\ga^N\setminus \{\pmb 0\}$.
We illustrate this by distinguishing two cases: $\ga\in G$ and $\ga\not \in G$.

If $\ga$ belongs to $G$,
by \cite[Lem.~3.2.1]{FJR1},
the condition $\vec n(\pmb x)=\pmb 0$  for $\pmb x\in{\CC_\ga^N}$
implies $\vec n(\pmb x)=\pmb 0$ for $\pmb x\in \CC^N$; hence
we have $\pmb x=\pmb 0$. In other words the hypersurface
$\{W_\ga=0\}$ inside $\CC_\ga^N\setminus \{\pmb 0\}$ is smooth.

On the other hand if $\ga\not \in G$, then
$\ga=(g_1 \la^{w_1},\dots,g_N \la ^{w_N})$
with $\la\not \in \pmmu_d$ and $(g_1,\dots, g_N)\in G$ (see \eqref{eq:barla} and \eqref{eq:intgroups}).
In this case $W_\ga$  vanishes identically on $\CC_\ga^N$. Indeed
suppose by way of contradiction
that ${x_1}^{m_1}\cdots x_q^{m_q}$
is a nonzero monomial of $W$ involving only $\ga$-fixed variables (\emph{i.e.} $g_1\la^{w_1}x_1=x_1$, \dots,
$g_q\la^{w_q}x_q=x_q$). Then $\la^d=1$ because we have
\begin{multline}\label{eq:restvanish}
x^{m_1}_1\cdots x_q^{m_q}=(g_1\la^{w_1}x_1)^{m_1}\cdots (g_q\la^{w_q}x_q)^{m_q}\\=
\la ^{w_1m_1+\cdots+ w_qm_q}((g_1x_1)^{m_1}\cdots (g_q x_q)^{m_q})=\la ^d(x_1^{m_1}\cdots x_q^{m_q}).
\end{multline}
A contradiction.

In this way,
a sector is
attached to each $\ga \in G\CC^\times$ and its
coarse space is always a quotient of a smooth variety
\begin{equation}
\label{eq:dichotomy}
\begin{cases}\{W_\ga =0\}_\ga/G\CC^{\times} \subset(\CC_\ga^N\setminus\{\pmb 0\})/G\CC^{\times}
&\text{if $\ga \in G$}.\\
\{W_\ga =0\}_{\ga}/G\CC^{\times}=(\CC_\ga^N\setminus\{\pmb 0\})/G\CC^{\times}
&\text{if $\ga \not \in G$;} \end{cases}
\end{equation}
\begin{rem}
   The second case of the above dichotomy
   corresponds to the situation where the intersection between
   $X_W$ and a twisted sector of
   the ambient space
   is not transverse. In fact,
   $X_W$ contains the twisted sector. This is the main difference between the Gorenstein and nonGorenstein cases,
   see Example \ref{exa:1nonG}.
   For a while, we considered it to be a major obstacle for the LG/CY correspondence.
   \end{rem}

The
action of $\ga$ on a fixed point $\pmb x\in \{W_\ga
=0\}_{\CC^N_\ga\setminus\{\pmb 0\}}$ on the tangent space
$T_{\pmb x}(\{W=0\})$ can be written (in a
suitable basis) as a diagonal matrix
$$\Diag(\exp(2\pi\cxi a^\gamma_1), \dots, \exp(2\pi\cxi a^\gamma_{N-1}))$$
for $a^\gamma_j\in [0,1[$.
Note that the matrix above is  $(N-1)\times(N-1)$
because $\{W=0\}$ is a smooth hypersurface in $\CC^N\setminus \{\pmb 0\}$.
We can read from the above matrix the so-called \emph{age shift}
\begin{equation}\label{eq:age}
a(\ga)=a\left(\Diag(\exp(2\pi\cxi a^\gamma_1), \dots, \exp(2\pi\cxi a^\gamma_{N-1}))\right)=\textsum_{l=1}^{N-1} a^\gamma_{l}.
\end{equation}
Note that here we regarded $\ga$ inside $GL(T_{\pmb x}(\{W=0\}),N-1)$, but in
our situation $\ga$ naturally operates
also on the affine space $\CC^N$;
we refer to Lemma \ref{lem:normal} in Section \ref{sect:main} for a formula expressing
the age $a_{\pmb x}(\ga)$ given above in terms
of the age of $\gamma$ as an element of $GL(\CC^N,N)$.

We finally define the bigraded Chen--Ruan cohomology as a direct sum of
ordinary cohomology groups of twisted sectors
\begin{equation}\label{eq:CRdefn}
H_{\CR}^{p,q} ([X_W/\wt G];\CC)=\bigoplus_{\gamma \in G\CC^\times}
H^{p-a(\ga),q-a(\ga)} (\{W=0\}_\ga/G\CC^{\times};\CC),
\end{equation}
where $\{W=0\}_\ga$ denotes the locus $\{\pmb x\in \{W=0\}_{\CC^N\setminus \{0\}}\mid \ga\pmb x=\pmb x\}$,
and the quotients appearing on the right hand side are quotient schemes and will be referred to
as sectors.
The total degree $\deg_{\CR}$ of a class $\al\in H_{\CR}^{p,q}
([X_W/\wt G];\CC)$ is $p+q$:
$$H_{\CR}^{d}([X_W/\wt G];\CC)=\bigoplus_{p+q=d}H^{p,q}_{\CR}([X_W/\wt G];\CC).$$
We do not discuss the Chen--Ruan orbifold product, because we only
regard $H_{\CR}$ as a bigraded vector space.

%the nondegenerate
%
%inner product is induced by
%the natural isomorphism between $\{W=0\}_\gamma$ and $\{W=0\}_{\gamma^{-1}}$
%$$H^{d_1}(\{W=0\}_\ga/G\CC^{\times};\CC)\otimes H^{d_2} (\{W=0\}_{\ga^{-1}}/G\CC^{\times};\CC)
%\longrightarrow \CC.$$
%Due to the age shift this yields indeed a duality between $H^{d}_{\CR}$ and $H^{2N-4-d}_{\CR}$,
%where $2N-4$ is the real dimension of $X_W$.

\subsection{The Landau--Ginzburg side}\label{sect:LGside}
On the Landau--Ginzburg side $W$ is regarded as a $G$-invariant
function $$W\colon \CC^N\to \CC,$$ and the fibre over the origin
is singular. We associate a nondegenerate bigraded vector space to
this singularity: the \emph{Fan--Jarvis--Ruan--Witten state
space}. It will be the counterpart on the Landau--Ginzburg side of
Chen--Ruan cohomology on the Calabi--Yau side.

For each $\gamma=(\exp(2\pi \cxi \Theta^\gamma_1),\dots,\exp(2\pi\cxi
\Theta^\gamma_N))\in G,$  with $\Theta_{\ga}^j\in [0,1[$; recall
the notations $\CC^N_{\gamma}$, $N_\gamma$, and $W_{\gamma}$ from
(\ref{eq:fixlocus}-\ref{eq:restfix}). The only critical point of
$W_\gamma$ is the origin (see \cite[Lem.~3.2.1]{FJR1}). Let
$\Hcal_{\gamma}$ be the $G$-invariant terms of the
middle-dimensional relative cohomology of $\CC^N_{\gamma}$
$$\Hcal_\ga=H^{N_\ga}(\CC^N_{\gamma}, W^{+\infty}_{\gamma};\CC)^{G},$$
where $W^{+\infty}=({Re} W_{\gamma})^{-1}\left(]M, +\infty[\right)$ for $M\gg0$.
The Fan--Jarvis--Ruan--Witten state space is
$$H_{\RW}(W,G;\CC)=\bigoplus_{\ga\in G} \Hcal_{\ga};$$
by analogy with Chen--Ruan cohomology, the summands will be often referred to as \emph{sectors}.
We point out a special sector: for $\ga=J_W$ the term $\Hcal_{\ga}$ is $1$-dimensional; indeed
$N_{J_W}=0$ and the relative cohomology has a single ($G$-invariant) generator $\mathsf 1_{J_W}$.
This is a good spot to introduce the so called \emph{Neveu--Schwarz} sectors:
\begin{defn}\label{defn:NS}
A sector $\Hcal_\ga$ is a \emph{Neveu--Schwarz} sector as soon as
$N_\ga$ vanishes\footnote{We refer to Example \ref{exa:homog} and Step 3 of the Proof of
Theorem \ref{thm:main} (Section \ref{sect:main}) for a geometric
interpretation of these sectors on the Calabi--Yau side.}. A
Neveu--Schwarz sector $\Hcal_{\ga}$ has  a single canonical
generator $\mathsf 1_{\ga}$. Following established practice we
call the remaining sectors \emph{Ramond} sectors (see
\cite{FJR1}).
\end{defn}

Using the Hodge decomposition of $\Hcal_\ga$
we define a bigraded decomposition of $H_{\RW}$.
As in Chen--Ruan cohomology, the age shift \eqref{eq:age} plays a role:
for example the total degree $d_{\RW}$ of the
terms $\Hcal_{\ga}$ is  equal to $N_\ga-2+2a(\ga)$
rather than the ordinary relative cohomology degree
$N_\ga$.
More precisely the decomposition of $\Hcal_{\ga}$ in terms of $\Hcal_{\ga}^{p,q}$ is as follows
\begin{align*}
&\Hcal_{\gamma}^{p,q}:= H^{p+1-a(\ga),q+1-a(\ga)}(\CC^N_{\gamma}, W^{+\infty}_{\gamma};\CC)^{G},\\
&\Hcal_{\ga}=\bigoplus_{p+q=N_\ga-2+2a(\ga)}\Hcal^{p,q}_\ga.
\end{align*}
The {state space} of FJRW theory is then equipped with a bigrading
    \begin{equation}\label{eq:FJRWdefn}
    H_{\RW}^{p,q}(W,G;\CC)=\bigoplus_{\gamma\in G} \Hcal^{p,q}_{\gamma}\end{equation}
and the total degree $\deg_{\RW}$ of a class in $H_{\RW}^{p,q}(W,G;\CC)$ is $p+q$;
note that, by construction, for any $\alpha \in \Hcal_\gamma$ and $\beta\in
\Hcal_{\gamma^{-1}}$ we have
\begin{equation*}\label{eq:degRel}
\deg_W(\alpha)+\deg_W(\beta) = 2N-4.
\end{equation*}
%Due to the canonical isomorphism between $\Hcal_{\ga}$ and $\Hcal_{\ga^{-1}}$ (see
%\cite{FJR1}),
%the state space $H_{\RW}(W,G)$ carries a natural Poincar\'e pairing
%    \begin{equation*}\label{eq:FJRWparing}
%    (\ ,\ )_{\RW}\colon \Hcal_{\gamma}\otimes \Hcal_{\gamma^{-1}}\rightarrow \CC.
%    \end{equation*}
%The pairing is nondegenerate and, due to the role of
%the age shift in the definition of $\deg_{\RW}$, yields a duality between $H_{\RW}^d$ and
%$H^{2N-4-d}_{\RW}$.

\begin{rem}\label{rem:CYanalogue} We make an observation which may be regarded as
the LG analogue of Remark \ref{rem:CY}.
The CY condition plays a crucial role here: the $\RW$-degree
of the canonical generator
$\mathsf 1_{J_W}$ of $\Hcal_{J_W}$ vanishes. We mention in
passing that, when the product is introduced, $\mathsf 1_{J_W}$ may be
regarded as a unit of $H_{\RW}(W,G;\CC)$ (see \cite{FJR1} and
\cite{Kr}).

Furthermore, in \cite{FJR1} the above structure is defined beyond
the case of the CY condition: it is important to notice that in
order to extend the structure together with the property
$\deg_{\RW}(\mathsf 1_{{J}})=1$ the authors involve  the charges
$q_1,\dots,q_N$ directly in the definition of the age shift (see
\cite[Defn.~3.2.3]{FJR1}).
\end{rem}

\subsection{The isomorphism}
The main theorem provides an
isomorphism between the Landau--Ginzburg side and the
Calabi--Yau side. As mentioned in the introduction this goes beyond the
expected correspondence for $G$ satisfying $J_W\in G\subseteq SL_W$ (see
Example \ref{exa:demo}, where $G\not \subseteq SL_W$).

\begin{thm} \label{thm:main}
Let $W$
be a nondegenerate
quasihomogeneous polynomial of degree $d$ in the variables $x_1,\dots,x_N$
whose charges add up to $1$ (CY condition).
Then, for any group $G$ of diagonal symmetries containing $J_W$
we have a bidegree-preserving isomorphism of vector spaces
$$H_{\CR}^{p,q}([X_W/\wt G];\CC)\cong H^{p,q}_{\RW}(W,G;\CC).$$
\end{thm}

For a scheme-theoretic counterpart of the above theorem we
should consider $G\subseteq SL_W$. Then we have the following statement.
\begin{cor}\label{cor:LG/CYscheme} Let $G$ be a subgroup of $SL_W$.
Assume that $X_W/\wt G$ admits a crepant resolution $Z$. Then,
we have $H^{p,q}(Z;\CC)\cong H^{p,q}_{\RW}(W,G;\CC).$
\end{cor}

See Section \ref{sect:main} for the proofs; we now discuss some examples.

\section{A first approach}\label{sect:warmup}

The following examples will provide a concrete
introduction to $\CR$ orbifold cohomology and the
$\RW$ state space. In each case we will establish by hand  the isomorphism
of Theorem \ref{thm:main} stated in the introduction.
This illustrates how certain sectors of the $\RW$ state space on the
Landau--Ginzburg side are interchanged
with cohomology classes on the Calabi--Yau side.
The exchange is nontrivial and
provides some early motivation for
the introduction of a bookkeeping device: the diagram introduced in
Section \ref{sect:main}. All the examples below will be examined
in Section \ref{sect:exa} using the diagram.

\begin{exa}[homogeneous polynomials]\label{exa:homog}
Theorem \ref{thm:main} is rather straightforward for a degree-$d$
hypersurface in $\PP^{d-1}$. Here $(w_1,\dots,w_d)$ is the $d$tuple
$(1,\dots,1)$ and the CY condition is automatically satisfied, $d=\sum_j w_j$.
This is the case of a cubic curve in $\PP^2$,
a K3 surface in $\PP^3$ (degree 4),
and a quintic three-fold in $\PP^4$.
The Lefschetz hyperplane theorem yields $N-1$ cohomology classes:
$\textsf 1\cap X_d,
\textsf h\cap X_d, \dots
\textsf h^{d-2}\cap X_d$ of bidegrees $(0,0), (1,1),\dots,(d-2,d-2)$.
The remaining classes, the cokernel of $H^\bullet(\PP^{d-1};\CC)\to H^\bullet(X_d;\CC)$,
are the primitive cohomology classes of degree $d-2$:
the $(p,q)$ primitive cohomology classes
can be identified with the ${J_W}$-invariant
$(p+1,q+1)$-classes of
$H^{d}(\CC^d,W^{+\infty}; \CC)$. For the cubic curve we have $(h^{1,0},h^{0,1})=(1,1)$,
for the K3 surface we have $(h^{2,0},h^{1,1},h^{0,2})=(1,20,1)$, and
for the quintic three-fold we have
$(h^{3,0},h^{2,1},h^{1,2},h^{0,3})=(1,101,101,1)$.
The Hodge ``diamond'' for the quintic polynomial
$(W,\langle {J_W}\rangle)$ on the
Calabi--Yau side (recall that $\wt{ \langle {J_W}\rangle}$
is the trivial group  $\langle {J_W}\rangle/ \langle {J_W}\rangle$) is \eqref{eq:quinticdiamond}.

If we switch to the Landau--Ginzburg side and we compute the
$\RW$ state space for $(W,\langle {J_W}\rangle)$, we
get $$H_{\RW}=\bigoplus_{i=0}^{d-1}\Hcal_{{J}^i}.$$
There are $d-1$ sectors, $\Hcal_{{J}^{i}}$ with $i\neq 0$,
for which $N_{{J}^i}$ vanishes: these are ${J_W}$-invariant relative cohomology classes
of bidegree $(0,0)$ in $H^{N_{{J}^i}}(\CC^{N_{{J}^i}},\varnothing;\CC)$. In other words
we have $d-1$ Neveu--Schwarz generators $\mathsf 1_{{J}}, \mathsf 1_{{J}^2}, \dots, \mathsf 1_{{J}^{d-1}}$
of $\RW$ bidegree $(0,0), (1,1), \dots, (d-2,d-2)$.
The sector $\Hcal_{{J}^0}$ is by definition
the ${J_W}$-invariant part of $H^{d}(\CC^d,W^{+\infty};\CC)$; therefore we get the same
Hodge diamond as on the Calabi--Yau side; \emph{i.e.} for the quintic three-fold we
get \eqref{eq:quinticdiamond}.

We can further test Theorem \ref{thm:main} by choosing a larger group $G\supsetneq \langle {J_W}\rangle$.
We will detail one of these calculations in Example \ref{exa:groupquot}.

There is only one observation that we wish to retain from this example:
\emph{the Neveu--Schwarz sectors on the LG side are interchanged
with the hyperplane sections on the CY side}. Note also that their degrees match.
\end{exa}
\begin{exa}[quasihomogeneous polynomials inside a Gorenstein $\PP(\pmb w)$]\label{exa:1inG}
Let us consider $W(x_1,x_2,x_3,x_4)=x_1^6+x_2^4+x_3^4+x_4^3$, which is quasihomogeneous
of degree $12$ in four variables of weight $2,3,3,4$.
On the Calabi--Yau side, we are interpreting this datum as a K3 surface $S$ inside the
Gorenstein weighted projective stack $\PP(2,3,3,4)$ (all weights divide the sum of the weights
$12$). We point out that the surface $S$ has only
two types of stack-theoretic points
with nontrivial stabilizers: the $3$ intersections of
$\{W=0\}$ with $\{x_2=x_3=0\}$, which have stabilizer $\pmmu_2$,
and the $4$ intersections of $\{W=0\}$ with
$\{x_1=x_4=0\}$, which have  stabilizer $\pmmu_3$.
These points contribute to
the twisted sectors: on the one hand
a point
$p$ with stabilizer $\pmmu_2$ yields the pair (point, automorphism)=$(p,1)$
in the ``untwisted'' sector $S_1$
and the pair $(p,-1)$ in the twisted sector $S_{-1}$, on the other hand
a point
$p$ with stabilizer $\pmmu_3$ yields $(p,1)$ in the ``untwisted'' sector $S_1$
and $(p,\xi_3)$ in the twisted sector $S_{\xi_3}$,
and $(p,\xi_3^2)$ in the twisted sector $S_{\xi_3^2}$.
In this way the ``twisted'' sectors ($S_{\ga}$ with $\ga\neq 1$)
consist of $4+4+3=11$ points.
It is straightforward to
see that all these points have age $1$:
therefore they contribute to an $11$-dimensional subspace of $H^{1,1}$
in $\CR$ orbifold cohomology. The remaining $\CR$ cohomology generators come from the
sector $S_1$, whose Hodge numbers are $(h^{2,0},h^{1,1},h^{0,2})=(1,9,1)$. Putting everything together,
we get the K3 surface Hodge diamond
\begin{equation}\begin{matrix}\label{eq:K3diamond}
& & &   &1&   &   \\
& & &0  & &0  &   \\
& &1&   &20&   &1  \\
& & &0  & &0  &   \\
\quad & & &   &1&   &&&.
\end{matrix}\end{equation}

On the LG side we compute the $\RW$ state space. There are
$12$ sectors
\[\begin{tabular}{c||c|c|c|c||c|c|}
   ${J_W}^h$& $x_1$ &$x_2$ &$x_3$ & $x_4$ &     $\deg_{\RW}$& $(h^{p,q}\mid p+q=\deg_{\RW})$\\
\hline&&&&&&\\
   ${J_W}^{0}$  &$0$      &$0$  &    $0$   & $0$    & 2    & $(h^{2,0},h^{1,1},h^{0,2})=(1,8,1)$ \\
   ${J_W}^{1}$  &$2$  &$3$  &$3$  & $4$ & 0    & $h^{0,0}=1$ \\
   ${J_W}^{2}$  &$4$  &$6$  &$6$  & $8$ & 2 & $h^{1,1}=1$\\
   ${J_W}^3$    &$6$  &$9$  &$9$  & $0$    &   &\\
   ${J_W}^4$  & $8$  &$0$     &$0$ & $4$ & 2 & $h^{1,1}=3$\\
   ${J_W}^5$  & $ 10$     &$3$ &$3$ & $8$ & 2& $h^{1,1}=1$\\
   ${J_W}^6$  & $ 0$     &$6$ &$6$ & $0$ & 2& $h^{1,1}=2$\\
   ${J_W}^7$  & $ 2$     &$9$ &$9$ & $4$ & 2& $h^{1,1}=1$\\
    ${J_W}^8$  & $4$  &$0$     &$0$ & $8$ & 2 & $h^{1,1}=3$\\
    ${J_W}^9$    &$6$  &$3$  &$3$  & $0$    &   &\\
    ${J_W}^{10}$  &$8$  &$6$  &$6$  & $4$ & 2 & $h^{1,1}=1$\\
    ${J_W}^{11}$  &$10$  &$9$  &$9$  & $8$ & 0    & $h^{2,2}=1$
\end{tabular}\ ,
\]
where the entry $m$ for a coordinate stands for
a coordinate $\exp(2\pi\cxi m)$
of the power of ${J_W}$ which we are considering. (We have put no entries
where there is no invariant element.)
Putting everything together we recover the same Hodge diamond \eqref{eq:K3diamond}.

We can test this further with the degree-$60$ three-fold
$\{x_1^{20}+x_2^6+x_3^5+x_4^4+x_5^3=0\}$ contained in $\PP(3,10,12,15,20)$.
We leave to the reader this interesting case, see Figure \ref{fig:bigCY} at the end of the paper.
The main point we wish to observe at this stage is that we find again the correspondence between
Neveu--Schwarz sectors and hyperplane generated cohomology classes. This is less
obvious than in the previous example because \emph{hyperplane generated classes
occur also in the twisted sector}:
for instance the sector $S_{-1}$ has $3$-dimensional and
corresponds to $\Hcal_{{J}^6}$ for the primitive part
and to one of the Neveu--Schwarz sectors for the nonprimitive part.
\end{exa}
\begin{exa}[a nonGorenstein ambient space $\PP(\pmb w)$]\label{exa:1nonG}
We now consider the polynomial $W=x_1^4x_2+x_2^3x_3+x_3^3x_4+x_4^3$
of degree $27$ and weights $5,7,6,9$.
On the CY side we have a K3 surface $S$ inside the nonGorenstein weighted projective stack $\PP(5,7,6,9)$.
The study of the special points whose stabilizer is nontrivial is rather subtle.
The ambient weighted projective stack has a point with stabilizer $\pmmu_9$ and a point with stabilizer $\pmmu_5$.
These two fixed loci behave differently with respect to $\{W=0\}$ and illustrate the dichotomy
\eqref{eq:dichotomy}: the first one $\{x_2=x_3=x_4=0\}$
is intersected transversely (\emph{i.e.} the intersection is empty because $\{x_2=x_3=x_4=0\}$ is a point),
the second one
$\{x_1=x_2=x_3=0\}$
is intersected nontransversely (\emph{i.e.} it is contained in $\{W=0\}$).
In Lemma \ref{lem:normal}
we show that this happens because the first stabilizer is an element of $\langle {J_W}\rangle $
whereas the second stabilizer is not.
This phenomenon is the crucial point of this example and may
be phrased as follows.

The stabilizers
$\pmmu_7$, $\pmmu_6$, and $\pmmu_5$ arise as subgroups of $\CC^\times$
generated by $\xi_7,\xi_6,\xi_5$
acting as $\la(x_1,\dots,x_4)=(\la^5x_1,\la^7x_2,\la^6x_3,\la^9x_4)$.
These elements are not contained in
the group generated by ${J_W}=(\xi_{27}^5,\xi_{27}^7,\xi_{27}^6, \xi_{27}^9)$. These
special group elements should be treated in a special way both on the CY side and the LG side.
This happens whenever the ambient space is not Gorenstein
and will require the study of
extra group elements (beyond $\langle {J_W}\rangle$)
(see Example \ref{exa:diagK3nonG} and
Figure \ref{fig:K3nonG}
 illustrating the present example).

We continue the computation, which yields the Hodge diagram for K3 surfaces
\eqref{eq:K3diamond}.
Indeed, the untwisted sector has one $(0,0)$-class, one $(2,2)$-class and
the following decomposition in degree two, $(h^{2,0},h^{1,1},h^{0,2})=(1,3,1)$. On the other hand
there are four
special points with stabilizers of order
$5,7,6,$ and $3$: namely $\{x_2=x_3=x_4=0\}$ (order 5),
$\{x_1=x_3=x_4=0\}$ (order 7), $\{x_1=x_2=x_4=0\}$ (order 6) and
$\{x_1=x_2=x_3^3+x_4^2=0\}$ (order 3).
These contribute to the twisted sectors with $(5-1)+(7-1)+(6-1)+(3-1)=17$
points representing $(1,1)$-classes due to the age shift (which is again $1$).
This matches \eqref{eq:K3diamond}.

%\begin{figure}%[h]
%% TeXgraph version 1.94 beta-7.5
%\begin{tikzpicture}%
%\useasboundingbox (-3,0)--(2.5,2.2);
%\pgfsetroundjoin%
%%objet1  (Cercle)
%\pgfxyline(-3,0)(2.2,0)
%\pgfcircle[fillstroke]{\pgfxy(-1,0)}{.05cm}
%\pgfcircle[fillstroke]{\pgfxy(0,0)}{.05cm}
%\pgfcircle[fillstroke]{\pgfxy(1,0)}{.05cm}
%\pgfcircle[fillstroke]{\pgfxy(-1,0.33)}{.05cm}
%\pgfcircle[fillstroke]{\pgfxy(-1,0.66)}{.05cm}
%\pgfcircle[fillstroke]{\pgfxy(1,0.40)}{.05cm}
%\pgfcircle[fillstroke]{\pgfxy(-1,1)}{.05cm}
%\pgfcircle[fillstroke]{\pgfxy(0,0.66)}{.05cm}
%\pgfcircle[fillstroke]{\pgfxy(0,1.33)}{.05cm}
%\pgfcircle[fillstroke]{\pgfxy(1,0.80)}{.05cm}
%\pgfcircle[fillstroke]{\pgfxy(-1,1.33)}{.05cm}
%\pgfcircle[fillstroke]{\pgfxy(-2,0.28)}{.05cm}
%\pgfcircle[fillstroke]{\pgfxy(-2,0.57)}{.05cm}
%\pgfcircle[fillstroke]{\pgfxy(-2,0.85)}{.05cm}
%\pgfcircle[fillstroke]{\pgfxy(-2,1.14)}{.05cm}
%\pgfcircle[fillstroke]{\pgfxy(-2,1.42)}{.05cm}
%\pgfcircle[fillstroke]{\pgfxy(-2,1.70)}{.05cm}
%\pgfcircle[fillstroke]{\pgfxy(-2,0)}{.05cm}
%\pgfcircle[fillstroke]{\pgfxy(-1,1.66)}{.05cm}
%\pgfcircle[fillstroke]{\pgfxy(1,1.20)}{.05cm}
%\pgfcircle[fillstroke]{\pgfxy(1,1.60)}{.05cm}
%%\pgfputat{\pgfxy(4,0.1)}{\pgftext[left,top]{\color{black}\small \text{age=$0$}}}\pgfstroke
%%\pgfputat{\pgfxy(4,0.6)}{\pgftext[left,top]{\color{black}\small \text{age=$1$}}}\pgfstroke
%%\pgfputat{\pgfxy(4,1.1)}{\pgftext[left,top]{\color{black}\small \text{age=$1$}}}\pgfstroke
%%\pgfputat{\pgfxy(4,1.6)}{\pgftext[left,top]{\color{black}\small \text{age=$3$}}}\pgfstroke
%\end{tikzpicture}%
%\caption{the sectors of $S$.\label{fig:sectorsK3}}
%\end{figure}

On the LG side we only can run a simple check for brevity. The CY side shows
$16$ sectors, as many as the elements of
$\pmmu_7 \cup \pmmu_6\cup \pmmu_3\cup \pmmu_5$, which contribute with
$18$ hyperplane sections (because the untwisted sector is two-dimensional
and yields $\mathsf 1, \mathsf h, \mathsf h^2$).
We find $20$ corresponding Neveu--Schwarz
sectors on the LG side: ${J_W}^h$ for $h$ prime to $\deg(W)=27$.
%%$$h=1,2,4,5,7,8,10,11,13,14,16,17,19,20,22,23,25,26.$$
\end{exa}
\begin{exa}[group quotients]\label{exa:groupquot}
We conclude this first study of the claim of Theorem \ref{thm:main}
with an example where $G \supsetneq \langle {J_W}\rangle$.
As the previous section already shows, a detailed analysis of the twisted sectors
on the CY side may be very delicate. Fortunately,
the theory of elliptic curves provides
a very well known and illuminating example.
We mention that
this provides an example where the Landau--Ginzburg/Calabi--Yau
correspondence holds beyond $SL(3, \CC)$.

Let $W(x_1,x_2,x_3)=x_1^2x_2+x_2^2x_3+x_3^3$ and set $G$ equal to the
maximal group $\Aut(W)$, which is cyclic of order $12$ and is generated by
the element $(\exp(2\pi\cxi 1/12), \exp(2\pi\cxi 10/12), \exp(2\pi\cxi 4/12))$.
The hypersurface defined by $W=0$
is a cubic curve in $\PP^2$.
The group $\wt G=G/\langle {J_W}\rangle $ is cyclic of order $4$ and
the action fixes the point represented by $e_0:=\{x_2=x_3=0\}$ (over
this coordinate subspace the polynomial $W$ vanishes).
We may regard $E=\{W=0\}$ as a genus-$1$ curve with a
marking $e_0\in E$:
an elliptic curve $(E,e_0)$.
Since there is only one elliptic curve
with automorphism group of order $4$ ($j$-invariant $1728$), we
know that $(E,e_0)$ is isomorphic to
$$(\CC/(\ZZ+\cxi \ZZ), [0]\in \CC)$$
and the automorphism may be
regarded as the complex multiplication by $\cxi$.
There are only three special orbits which do not consist of four distinct points:
the one-point orbit $\{e_0=[0]\}$ (with stabilizer $\wt G$),
the one-point orbit $\{1/2+\cxi/2\}$ (with stabilizer $\wt G$),
and the two-point orbit containing $1/2$ and $\cxi/2$ (both with stabilizer of order $2$).
Therefore the stack-theoretic quotient $[E/\wt G]$
has only three special (\emph{i.e.} nonrepresentable) points
with stabilizers of order $m_0=4$, $m_1=4$, and $m_2=2$ (the
coarse space is actually a projective line $E/\wt G\cong \PP^1$).
It is now easy to visualize the sectors: apart from the ``untwisted'' sector,
we find $\sum_i (m_i-1)=7$ ``twisted'' sectors
corresponding to points paired with their nontrivial automorphism.
We expect a $9$-dimensional $\CR$ cohomology
vector space $H^{\bullet}_{\CR}([E/\wt G];\CC)$
with a $2$-dimensional contribution from the ``untwisted'' sector ($H^\bullet(\PP^1)\cong \mathsf 1\CC\oplus \mathsf h\CC$)
and seven twisted $1$-dimensional contributions mentioned above (graded by twice the age).
The picture is illustrated in Figure \ref{fig:sectorsorbicurve}, where the Hodge numbers are also listed.
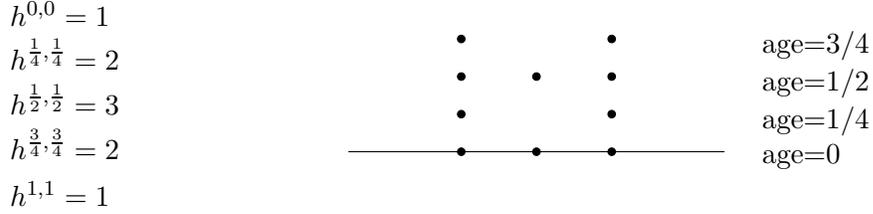
\begin{figure}
% TeXgraph version 1.94 beta-7.5
\begin{tikzpicture}%
\useasboundingbox (-5.5,-1)--(2.5,2.5);
\pgfsetroundjoin%
%objet1  (Cercle)
\pgfxyline(-2.5,0)(2.5,0)
\pgfputat{\pgfxy(-7,-0.4)}{\pgftext[left,top]{\color{black}\small \text{$h^{1,1}=1$}}}\pgfstroke
\pgfputat{\pgfxy(-7,0.3)}{\pgftext[left,top]{\color{black}\small \text{$h^{\frac34,\frac34}=2$}}}\pgfstroke
\pgfputat{\pgfxy(-7,.9)}{\pgftext[left,top]{\color{black}\small \text{$h^{\frac12,\frac12}=3$}}}\pgfstroke
\pgfputat{\pgfxy(-7,1.5)}{\pgftext[left,top]{\color{black}\small \text{$h^{\frac14,\frac14}=2$}}}\pgfstroke
\pgfputat{\pgfxy(-7,2.0)}{\pgftext[left,top]{\color{black}\small \text{$h^{0,0}=1$}}}\pgfstroke
\pgfcircle[fillstroke]{\pgfxy(-1,0)}{.05cm}
\pgfcircle[fillstroke]{\pgfxy(0,0)}{.05cm}
\pgfcircle[fillstroke]{\pgfxy(1,0)}{.05cm}
\pgfcircle[fillstroke]{\pgfxy(-1,0.5)}{.05cm}
\pgfcircle[fillstroke]{\pgfxy(1,0.5)}{.05cm}
\pgfcircle[fillstroke]{\pgfxy(-1,1)}{.05cm}
\pgfcircle[fillstroke]{\pgfxy(0,1)}{.05cm}
\pgfcircle[fillstroke]{\pgfxy(1,1)}{.05cm}
\pgfcircle[fillstroke]{\pgfxy(-1,1.5)}{.05cm}
\pgfcircle[fillstroke]{\pgfxy(1,1.5)}{.05cm}
\pgfputat{\pgfxy(3,0.1)}{\pgftext[left,top]{\color{black}\small \text{age=$0$}}}\pgfstroke
\pgfputat{\pgfxy(3,0.6)}{\pgftext[left,top]{\color{black}\small \text{age=$1/4$}}}\pgfstroke
\pgfputat{\pgfxy(3,1.1)}{\pgftext[left,top]{\color{black}\small \text{age=$1/2$}}}\pgfstroke
\pgfputat{\pgfxy(3,1.6)}{\pgftext[left,top]{\color{black}\small \text{age=$3/4$}}}\pgfstroke
\end{tikzpicture}%
\caption{the sectors of $[E/\wt G]$.\label{fig:sectorsorbicurve}}
\end{figure}
\end{exa}

We finally check that the above computation matches the LG side. By $\ga$, we denote
the order-$12$ generator of $G$.
\[\begin{tabular}{c||c|c|c||c|c|}
   $\ga^h$& $x_1$ &$x_2$ &$x_3$  &     $\deg_{\RW}$& $(h^{p,q}\mid p+q=\deg_{\RW})$\\
\hline&&&&&\\
   $\ga^{0}$  &$0$      &$0$  &    $0$      &  &  \\
   $\ga^{1}$  &$1$  &$10$  &$4$             & 1/2 & $h^{1/4,1/4}=1$ \\
   $\ga^{2}$  &$2$  &$8$  &$8$              & 1 & $h^{1/2,1/2}=1$\\
   $\ga^3$    &$3$  &$6$  &$0$              &   &\\
   $\ga^4$  & $4$  &$4$     &$4$            & 0 & $h^{0,0}=1$\\
   $\ga^5$  & $5$     &$2$ &$8$             & 1/2 & $h^{1/4,1/4}=1$\\
   $\ga^6$  & $6$     &$0$ &$0$             & 1 & $h^{1/2,1/2}=1$\\
   $\ga^7$  & $7$     &$10$ &$4$            & 3/2 & $h^{3/4,3/4}=1$\\
    $\ga^8$  & $8$  &$8$     &$8$           & 2 & $h^{1,1}=1$\\
    $\ga^9$    &$9$  &$6$  &$0$             &   &\\
    $\ga^{10}$  &$10$  &$4$  &$4$           & 1 & $h^{1/2,1/2}=1$\\
    $\ga^{11}$  &$11$  &$2$  &$8$           & 1 & $h^{3/4,3/4}=1$
\end{tabular}\
\]
Once again we put no
entries where there is no invariant element. The Hodge numbers
match those listed in Figure \ref{fig:sectorsorbicurve}.

\section{Proof of the main result: a combinatorial model}\label{sect:main}
The proof is structured in five steps as follows. On the Calabi--Yau side,
we further detail the decomposition of the $\CR$ cohomology (Step 1).
Then, we do the same for the $\RW$ state space on
the Landau--Ginzburg side (Step 2).
We provide a diagram
which schematizes and assembles into one picture
the sectors on the two sides (Step 3).
We prove a lemma which allows
us to read off $\deg_{\CR}$ and $\deg_{\RW}$ on the diagram (Step 4).
We establish an isomorphism using the combinatorial model (Step 5).

\medskip

\paragraph{Step 1: \emph{Calabi--Yau side}.} Consider the decomposition
 \eqref{eq:CRdefn} of $H_{\CR}$ as a sum over $G\CC^\times$. The complex dimension
 of $H_{\CR}$ is finite although this is not evident from \eqref{eq:CRdefn}.
 Indeed, we can decompose $G\CC^\times$ modulo $\CC^{\times}$
 into $M=\abs{G}/d$ cosets. Let us choose $M$ distinct
 cosets $g^{(1)}\CC^\times,\dots, g^{(M)}\CC^\times$ so that
 $g^{(1)},\dots, g^{(M)}\in G$ and
 the set $\sqcup_{i=1}^M g^{(i)}\CC^\times$ equals the set $G\CC^\times$. Now,
 we describe the direct sum
 \begin{equation}\label{eq:coset}
\bigoplus_{\ga\in g\CC^\times}H^{\bullet}(\{W_\ga=0\}_\ga/G\CC^\times ;\CC)
\end{equation}
 where $g$ is any of the elements $\{g^{(1)},\dots,g^{(M)}\}$.
 By construction $H_{\CR}$ is the direct sum of the expressions above for
 $g$ ranging over $\{g^{(1)},\dots, g^{(M)}\}$.

 Now we exhibit a finite number of terms of $g\CC^\times$, outside which
 the summand of \eqref{eq:coset} vanishes.
 Regard an element $g\in G$ as an $N$-tuple of elements of $\CC^\times$,
 $$g=(g_j)_{j=1}^N.$$
 Notice that specifying $\ga$ in $g\CC^\times$ is equivalent to
 choosing $\la\in \CC^\times$ so that $\ga=g\bar \la =(g_j)_{j=1}^N(\la^{w_j})_{j=1}^N$.
 Since $g\bar \la$ acts by multiplication on the coordinates, the fixed locus is nonempty
 if and only if $\la$ is contained in the finite set $\bigcup_{j=1}^N \{\la\mid \la^{-w_j}=g_j\}$.
 In this way \eqref{eq:coset} can be rewritten as a direct sum of a finite number of
 finite dimensional vector spaces
 \begin{equation}\label{eq:dottedrays}
\bigoplus_{\la \in \ \bigcup_{j=1}^N \left\{\la\mid \la^{-w_j}=g_j\right\}}
H^{\bullet}(\{W_{g\bar \la}=0\}_{g\bar \la}/G\CC^\times ;\CC),
\end{equation}
 where the notation $\bar \la$ of \eqref{eq:barla} has been used.

 The quotient scheme $\{W_{g\bar \la}=0\}_{g\bar \la}/G\CC^\times$ may be regarded
 as the quotient scheme by $G\CC^\times/\CC^\times  =\wt G$
 of the hypersurface $\{W_{g\bar\la}=0\}$ inside the weighted projective
 space $\PP(\pmb w_\la)$ where $\pmb w_\la$ is the multi-index
 \begin{equation}\label{eq:multiindexlambda}
 \pmb w_{\la}=\{w_j\mid  \la^{-w_j}=g_j\}.
 \end{equation}
 In this way we have
 $$H^{\bullet}(\{W_{g\bar \la}=0\}_{g\bar \la}/G\CC^\times ;\CC)=
 H^{\bullet}(\{W_{g\bar \la}=0\}_{\PP(\pmb w_\la)};\CC)^{\wt G}.$$
 Notice that the number of entries of $\pmb w_\la$ equals $N_{g\bar\la}$.

 The cohomology $H^\bullet$ of a hypersurface $S$ inside a
 weighted projective stack splits into  two summands.
 The first summand
 is generated by the self-intersections
 of the hyperplane sections: $\mathsf 1_S$, $\mathsf h\cap S$, $\mathsf h^2\cap S$, $\dots$
 In the case of $\{W_{g\bar \la}=0\}_{\PP(\pmb w_\la)}$,
 this summand of $H^\bullet(\{W_{g\bar \la}=0\}_{\PP(\pmb w_\la)};\CC)$ is
 $(N_{g\bar\la}-1)$-dimensional. We point out that all these terms are $\wt G$-invariant.
 The second summand is the primitive cohomology and is concentrated in
 degree $\delta=\dim_{\CC}(S)$ (if $\dim_\CC(S)$ is odd this summand is the entire
 cohomology group $H^\delta(\cdot \ ;\CC)$, otherwise the rank of this summand
 equals the Betti number $b_\delta=\dim H^\delta$ minus 1).
 By the theory of the Milnor fibre
 \cite{St} \cite{Do} \cite{Di} we may express the primitive cohomology as
 $H^{N_\ga}(\CC^N_{\gamma}, W^{+\infty}_{\gamma};\CC)^{\langle {J_W} \rangle}.$
This happens because the ${J_W}$-action  is the monodromy
action on the Milnor fibre of
$$W_{\ga}\colon\CC^{N_\ga}\to\CC.$$
In this way the $\wt G$-invariant part of
the primitive cohomology of the hypersurface $\{W_{g\bar \la}=0\}_{\PP(\pmb w_\la)}$
is isomorphic to $H^{N_\ga}(\CC^N_{\gamma}, W^{+\infty}_{\gamma};\CC)^G$ (the isomorphism identifies $(p,q)$-classes
in $H^{N_\ga}(\{W_{g\bar \la}=0\}_{g\bar \la}/G\CC^\times)$ with $(p+1,q+1)$-classes
in $H^{N_\ga}(\CC^N_{\gamma}, W^{+\infty}_{\gamma};\CC)^G$).
In this way the group $H^{\bullet}(\{W_{g\bar \la}=0\}_{g\bar \la}/G\CC^\times ;\CC)$ can
be decomposed as
\begin{equation}\label{eq:Lefschetzdecomp}
H^{N_{g\bar \la}}\left(\CC^N_{g\bar \la}, W^{+\infty}_{g\bar \la};\CC\right)^G\oplus \bigoplus_{i=0}^{N_{g\bar \la}-2} \left [\mathsf h^i\cap\{W_{g\bar \la}=0\}_{\PP(\pmb w_\la)}\right]\CC.
\end{equation}
\begin{rem}\label{rem:Hdgfilt}
The summands on the right hand side contain $(i,i)$-classes corresponding
to cohomology classes in $H^{2i}(\{W_{g\bar \la}=0\}_{g\bar \la}/G\CC^\times)$;
whereas the first summand consists of $(p+1,q+1)$-classes (with $p,q\ge 0$)
of degree $N_\gamma$ which represent $(p,q)$-classes in the primitive cohomology
of $\{W_{g\bar \la}=0\}_{g\bar \la}/G\CC^\times$.\end{rem}
By summing the above expression over all $\la \in \ \bigcup_{j=1}^N \left\{\la\mid \la^{-w_j}=g_j\right\}$
we  get the entire finite-dimensional
contribution to $H^{\bullet}_{\CR}$ coming from the coset $g\CC^\times$.

\medskip

\paragraph{Step 2: \emph{Landau--Ginzburg side}.} We analyze
the $\RW$ state space in a similar way:
$$    H_{\RW}^{\bullet}(W,G;\CC)=\bigoplus_{\gamma\in G} \Hcal_{\gamma}.$$
For $J=J_W$, we decompose $G$ into $M=\abs{G}/d$ distinct cosets
$g^{(1)}\langle {J}\rangle$,\dots, $g^{(M)}\langle {J}\rangle$ (we choose
the same $g^{(1)}, \dots, g^{(M)}$ as in the previous step).
Therefore the $\RW$ state space is a direct sum of the terms
$$\bigoplus_{i=0}^{d-1} H^{N_{g{J}^i}}\left(\CC^N_{g{J}^i}, W^{+\infty}_{g{J}^i};
\CC\right)^G$$
for $g$ ranging in $\{g^{(1)},\dots, g^{(M)}\}$ (we are just making the definition of $\Hcal_{g{J}^i}$ explicit).

Write $g=(g_j)_{j=1}^N$ as usual. We point out that
if $\xi_d^i$ does not belong to $\bigcup_{j=1}^N \left\{\la\mid \la^{-w_j}=g_j\right\}$, then
$N_{g{J}^i}=0$. In other words   $\Hcal_{g{J}^i}$ is of Neveu--Schwarz type.
We finally express the
entire contribution to $H^{\bullet}_{\RW}$ coming from the coset $g\langle {J}\rangle$:
$$\bigoplus_{\la \in \ \pmmu_d\cap \bigcup_{j=1}^N \left\{\la\mid \la^{-w_j}=g_j\right\}}
H^{N_{g\bar \la}}\left(\CC^N_{g\bar \la}, W^{+\infty}_{g\bar \la};\CC\right)^G\oplus\!\!\!\!\!\!\!\!\!\! \bigoplus_{\la \in \ \pmmu_d\setminus \bigcup_{j=1}^N \left\{\la\mid \la^{-w_j}=g_j\right\}}
\pmb 1_{g\bar \la}\CC,$$
where we used the notation \eqref{eq:barla},
and we identified the terms of $\langle {J}\rangle$
as $\bar \la$ for $\la \in \pmmu_d$ (\emph{e.g.} ${J}=\bar {\xi_d}$).

\medskip

\paragraph{Step 3: \emph{the diagram}.}
In the previous two steps we split the state spaces into $M$ summands
corresponding to
a set of $M$ elements  $g^{(1)},\dots,g^{(M)}$ in $G$.
Each summand is efficiently represented by a diagram, which may be regarded as
a generalization of
Boissi\`ere, Mann, and Perroni's model \cite{BMP}.

Again, let us choose one of the above elements $g^{(1)},\dots, g^{(M)}$
and denote it by $g$;
we describe the corresponding diagram. It consists of halflines (rays)
stemming from the origin in the complex plane and points lying on
them (dots). The dots will correspond to
(sets of) generators in $\CR$ cohomology, whereas
the rays will represent sectors of the $\RW$ state space.
Draw a \emph{ray}
$$\{\rho\nu \in \CC\mid \rho\in \RR^+\}\subset \CC$$
whenever we have
\begin{equation}\label{eq:setrays}
\nu \in
\pmmu_d\cup \bigcup_{j=1}^N \left\{\al\in \CC\mid \al^{w_j}=g_j\right\}.
\end{equation}
Mark a \emph{dot}
$$j \nu\in \CC$$
whenever $\nu^{w_j}=g_j$ for some $j$; in other words, whenever $\nu$ and $j$
satisfy
$$\nu \in \left\{\al\in \CC\mid \al^{w_j}=g_j\right\}.$$
Mark further dots
$$(N+1) \nu$$
whenever
$$\nu \in
\left( \bigcup_{j=1}^N \left\{\al\in \CC\mid \al^{w_j}=g_j\right\}\right)\setminus \pmmu_d.$$
For a nontrivial but low-dimensional example we refer the reader to
Figure \ref{fig:K3nonG}
where the diagram is drawn for
the above-mentioned K3 surface $\{x_1^4x_2+x_2^3x_3+x_3^3x_4+x_4^3\}\subset \PP(5,7,6,9)$.

This model  can be related to the sectors of the two
$\CR$ and $\RW$ spaces.
The coset determined by $h$ with $h=1\in G$ is the case treated in
\cite{BMP} and,  for the sake of clarity,  we discuss it first.
This corresponds to assuming $G=\langle {J}\rangle$ and looking at
the hypersurface $\{W=0\}\subset \PP(\pmb w)$ (if $G=\langle {J}\rangle$, then
$\wt G=1$).
Since $g_j=1$ for all $j$, following \eqref{eq:setrays}, we find that the rays correspond to
the elements of
$\pmmu_d\cup \pmmu_{w_1}\cup \dots\cup \pmmu_{w_N}.$
The rays that carry some dots
are in one-to-one correspondence with the sectors associated to
 the hypersurface $\{W=0\}$ inside
$\PP(w_1,\dots,w_N)$. If we write a  ray as
$\{\rho\nu\mid \rho\in \RR^+\}$ with $\abs{\nu}=1$ then the corresponding sector is
the hypersurface $\{W_\la=0\}_{\PP(\pmb w_\la)}$ for $\la=\nu^{-1}$.
Simply by unraveling the definitions, the authors of \cite{BMP}
make the following useful observation: a ray carries as many
dots as the quasihomogeneous coordinates of the corresponding
weighted projective subspace $\PP(\pmb w_\la)$. Building upon this,
one can derive a combinatorial model for the cohomology of the sector
$S=\{W_\la=0\}_{\PP(\pmb w_\la)}$: namely, we let the first $N_{\la}-1$ dots represent
the hyperplane sections $\mathsf 1_S$, $\mathsf h\cap S$, $\mathsf h^2\cap S$, $\dots$,
$\mathsf h^{N_\la-2}\cap S$ and
the $N_\la$th dots represent the primitive cohomology. In this way all the dots are attached to
a summand of the $\CR$ cohomology of $X_W$. On the Landau--Ginzburg side, we can use the diagram as follows:
the rays with
angular coordinate $2\pi l/d$ can be associated to the
summand $\Hcal_{{J}^{-l}}$ of the $\RW$ state space of $(W,G=\langle {J}\rangle)$.
The number of dots on one of these rays corresponds to the index $N_{{J}^{-l}}$.

The general procedure for a coset represented by $h$ is as follows.
Similarly to the case $h=1$, the rays whose angular coordinate is $2\pi l/d$
represent the sector of the $\RW$ space $\Hcal_{h{J}^{-l}}$.
We point out that, by construction, a sector is of Neveu--Schwarz type if and only if
it is \emph{empty}; \emph{i.e.} it does not carry any dot.
The dots always lie on some ray by construction: consider
the dot $m\nu$ (with $m\in \NN$ and $\mu\in \{z\mid \abs{z}=1\}$)
lying  on the ray $\{\rho\nu\mid \rho\in \RR^+\}$. We say that it is an
\emph{extremal dot} if there is no other dot with higher polar coordinate and is
an \emph{internal dot}
otherwise.
An extremal dot $m\nu$ corresponds to
the primitive cohomology of
$H^{\bullet}(\{W_{g\bar \la}=0\}_{g\bar \la}/G\CC^\times ;\CC)$
for $\la=\nu^{-1}$.
The internal dots $m_1\nu, m_2\nu, m_3\nu, \dots$ lying on $\{\rho \nu\mid \rho\in \RR\}$
can be ordered with respect to their polar coordinates and represent
hyperplane sections in Chen--Ruan cohomology
of the sector $\{W_{g\bar \la}=0\}_{g\bar \la}/G\CC^\times$ for $\la=\nu^{-1}$:
the first dot corresponds to
the fundamental class  of
$\{W_{g\bar \la}=0\}_{g\bar \la}/G\CC^\times$, the next corresponds
to the intersection with $\mathsf h$, and so on.

We refer to Example \ref{exa:demo} for a simple, and
nevertheless interesting, demonstration of
the above procedure (we wrote it in such a way that
the reader can skip directly there for a
detailed description of the diagram attached to
a coset).

Now, we define two functions $D$ and $R$
on the union of the sets of rays and of dots. They essentially count dots and rays and they
can be efficiently used in order to express the
quantities $\deg_{\CR}$ and $\deg_{\RW}$ for the corresponding classes.
Notice that dots and rays are naturally ordered:
the rays can be arranged according to the angular coordinate ranging over $[0,1[$
whereas the dots can be arranged in
lexicographic order $\preccurlyeq$ (recall that for $\vartheta,\vartheta'\in [0,1[$
we write $\rho\exp(2\pi\cxi\vartheta)\preccurlyeq \rho'\exp(2\pi\cxi\vartheta')$ if and only if
we have $\vartheta \le \vartheta'$ or, for $\vartheta=\vartheta'$, we have $\rho\le \rho'$).
We can actually order the set given by the union of dots and rays:
for this, we require
that a ray precedes all dots lying on it and on the following rays
(to this effect a ray $\{\rho \nu\mid \rho\in \RR^+\}$
may be treated as the point $(1/2)\nu$ and arranged according to
$\preccurlyeq$). Now we define the functions $R$ and $D$. The function
$R$ is naturally defined on all rays and takes values in the natural numbers ranging from $0$
to the size of the set
$\bigcup_{j=1}^N \left\{\al\in \CC\mid \al^{w_j}=g_j\right\}$
minus one. It is defined by simply counting the rays in
the sense of the angular coordinate (\emph{i.e.} anticlockwise). The function $D$ is naturally
defined on the set of dots and takes values in the natural numbers ranging from $0$ to the number
of dots minus $1$. It is defined by counting the dots in lexicographic order. We may naturally
extend the function $D$ to the set of rays: simply assign to a ray the value $D$ of the first preceding
dot (if the ray precedes all dots we set $D=-1$). We naturally extend $R$ to
the set of dots: a
dot takes the value $R$
assigned to the ray on which it lies.

\begin{rem}\label{eq:usingCYcond}
The functions $R$ and $D$ range over the same finite set of
numbers. This happens because the number of rays is clearly $d$
plus the number of elements of $\bigcup_j \pmmu_{w_j}\setminus
\pmmu_d.$ On the other hand the number of dots can be computed as
follows. The number of dots $j \nu$ with $\abs {j\nu}\le N$ is
$\sum_j w_j$ because each equation $\nu^{w_j}=g_j$  has $w_j$
solutions. The remaining dots are precisely as many as the
elements of $\bigcup_j \pmmu_{w_j}\setminus \pmmu_d$ by
construction. The two counts match by the CY condition: $d=\sum_j
w_j$.
\end{rem}

\medskip

\paragraph{Step 4: \emph{the degrees $\deg_{\CR}$ and $\deg_{\RW}$}.}
Let $\pmb x\in \CC^N$ be a point in $$({\CC}_{\ga}^N\setminus \{\pmb 0\})\cap \{W=0\}=\{W_\ga=0\}_\ga.$$
(\emph{i.e.} $\ga\pmb x=\pmb x$ and $W(\pmb x)=0$).
By \eqref{eq:dichotomy}, if $\ga\not \in G$, then the intersection is
\emph{not transversal}
and $\CC_\ga^N$ lies inside $\{W=0\}$; otherwise, if $\ga\in G$,
the intersection is transversal and the intersection locus is
again a smooth variety. Indeed, one can see  directly that if $\ga \in G$
the normal vector $\vec n(\pmb x)$ to $\pmb x \in \{W=0\}$ lies in $\CC_\ga^N$: hence the
whole
line $$\{\pmb y=\pmb x+\rho\vec{n}(\pmb x)\in \CC^N\mid \rho\in \RR\}$$
is fixed (lies inside $\CC^N_{\ga}$).

The explicit argument is as follows: let us
arrange the coordinates so that $x_1,\dots,x_q$ are all the
$\gamma$-fixed  coordinates: \emph{i.e.} if
$\ga=(g_1,\dots,g_N)$ we have $g_1= \dots =g_q=1$. Then,
for any $j>q$ we have $g_j\neq 1$. We conclude that $\partial_j W(\pmb x)=0$.
This happens because $\pmb x\in \CC_\ga^N$
is of the form $\pmb x=(x_1,\dots,x_q,0,\dots,0)$ and $\partial_jW(\pmb x)\neq 0$ only if
there is a monomial of $W$ of the form $x_1^{m_1}\cdots x_q^{m_q}x_j$, which contradicts
$g_j\neq 1$ because
$$x_1^{m_1}\cdots x_q^{m_q}x_j=(g_1x_1)^{m_1}\cdots (g_q x_q)^{m_q}(g_jx_j)=g_j(x_1^{m_1}\cdots x_q^{m_q}x_j).$$

In the case $\ga\not \in G$ we know that
the normal line passing through $\pmb x$ with vector $\vec n(\pmb x)$ has only one fixed point:
$\pmb x$. The following lemma describes this action precisely and
embodies the previous observation that $\ga$ acts trivially on $\pmb x$ for $\ga\in G$.
\begin{lem}\label{lem:normal}
For any $\ga=g\bar\la \in G\CC^{\times}$, let
$\pmb x\in \CC^N\setminus \{\bf 0\}$ be a point of the hypersurface
$\{W=0\}$, which is fixed by $\ga$; \emph{i.e.} $\pmb x$ belongs to
$({\CC}_{g\bar \la }^N\setminus \{{\bf0} \})\cap \{W=0\}.$
Then $g\bar\la$ acts on the normal line
$\{\textbf y=\pmb x+\rho\vec{n}(\pmb x)\in \CC^N\mid \rho\in \RR\}$
by multiplication by $\la^d$ as follows
$$g\bar \la \left(\pmb x+\rho\vec{n}(\pmb x)\right)=\pmb x+\la^d\rho\vec{n}(\pmb x).$$

In particular,
the age $\al$ of $g\bar \la$ in $GL(\CC,N)$ and the age $a_{\pmb x}(\ga)$ of $g\bar\la$
acting on the $(N-1)$-dimensional tangent space $T_{\pmb x}(\{W=0\})$ are related as follows:
$$a_{\pmb x}(g\bar \la)=\al -  \langle sd \rangle\qquad \text{if $\la=\exp(2\pi\cxi s)$ and $s\in [0,1[$},$$
where $\langle sd\rangle$ denotes the fractional part of $sd$  (\emph{i.e.} $sd-\lfloor sd\rfloor$).

As a consequence, on the diagram attached to
$g=(g_1,\dots, g_N)\in G$,
the degree $\deg_{\RW}$ of a class represented by an empty ray
and the degree $\deg_{\CR}$ of a class represented by an internal dot
can be expressed as
$$2\left(\sum_{j=1}^N s_j+D-R\right),$$
where $g_j=\exp((2\pi\cxi s_j))$ with $s_j\in [0,1[$.
\end{lem}
\begin{proof}
The first part is well known: the normal bundle to the hypersurface is
a $\CC^\times$-linearized line bundle
 $\Ocal(d)$ with
character $\la\mapsto \la ^d$.
We detail the argument by choosing
a nonvanishing coordinate $\partial_{j_0}W(\pmb x)$
of $\vec n(\pmb x)$ and by proving that multiplying it by
$g_{j_0}\la^{w_{j_0}}$ is the same as rescaling it by $\la^d$.
To begin with, notice that the fact that $\partial_{j_0}W(\pmb x)$
does not vanish guarantees the existence of a monomial of $W$
with exponents ${m_1},\dots, m_N$ only involving
the ${j_0}$th coordinate and coordinates for which
$g_j\la^{w_j}=1$. In other words, for $j\neq {j_0}$ we have
$(g_j\la^{w_j})^{m_j}=1$, because either $m_j$ vanishes or
$g_j\la^{w_j}$ equals $1$.

Then there are two possibilities.
First, if $g_{j_0}\la^{w_{j_0}}=1$, then $\la^d=1$,
$$\la^d=\la^{m_1w_1+\cdots+m_Nw_N}=g_1^{m_1}\la^{m_1w_1}\cdots g_N^{m_N}\la^{m_Nw_N}=(g_{j_0}\la^{w_{j_0}})^{m_{j_0}}=1.$$
Otherwise $g_{j_0}\la^{w_{j_0}}\neq 1$ and the $x_{j_0}$ coordinate is not $\ga$-fixed.
In this case $\partial_{j_0} W(\pmb x)\neq 0$ implies that
$m_{j_0}$ is necessarily equal to $1$: we have
$$g_{j_0} \la^{w_{j_0}}=g_{j_0} \la^{w_{j_0}}\prod_{j\neq {j_0}}(g_j\la^{w_j})^{m_j} = \prod_j g_j^{m_j} \prod_j \la^{m_j w_j}=\la^d.$$
This completes the proof of the first part of the claim.

The formula immediately implies the expression for $a_{\pmb x}(g\bar \la)$ in
terms of $\al$ and $\la$ in the statement.
Indeed, we make that expression more explicit
by assuming that $g$ equals $(\exp((2\pi\cxi s_j)))_{j=1}^N$ and
by writing $\la$ as $\exp(-2\pi\cxi t)$. Then we have
\begin{align*}
a_{\pmb x}(g\bar \la)&= \sum_{j=1}^N \langle s_j - tw_j\rangle - \langle-td\rangle\\
&=\sum_{j=1}^N (s_j - tw_j) - \sum_{j=1}^N \lfloor s_j - tw_j \rfloor - (-td) -(-\lfloor -td\rfloor)\\
&=\sum_{j=1}^N s_j +\left(- \sum_{j=1}^N \lfloor s_j - tw_j \rfloor\right) -\left(-\lfloor -td\rfloor\right),
\end{align*}
where the CY condition has been used in the last equality.
The last
part of the statement follows from relating the last two
summands to the function
$D-R$ evaluated on an empty ray and internal dots.

The functions $D$ and $R$ introduced above have particularly convenient
properties, which will be evident in the next step; however, in order to match the above
expression we need to define two slightly different functions $\wt D$ and $\wt R$.
The functions $\wt D$ and $\wt R$ only count (and are defined on) a special kind of dots and rays:
the rays are those  with angular coordinate within
$(2\pi/d)\NN$ and the dots are those whose polar coordinate  is (strictly) smaller than
$N+1$ (\emph{i.e.} $\abs{\cdot \ }\le N$).
The union of these dots and rays is naturally ordered by the
lexicographic order $\preccurlyeq$ and the prescription
that a ray precedes all dots lying on it and on the following rays.
The function
$\wt R$ is naturally defined on the considered rays by the angular coordinate times $d/2\pi$
and takes values in $\{0,1,\dots, d-1\}$.
The definition extends immediately to dots lying on the above-mentioned rays
and also
to a dot which does not lie on the considered rays: we assign
to it the value $\wt R$ of the next ray (and we assign $d$ if there is no next ray).
On the other hand, the function $\wt D$ is defined by counting in lexicographic order
the dots with $\abs{\cdot }\le N$.
Again, we may naturally
extend the function $\wt D$ to the set of rays: simply assign to a ray the value $\wt D$ of the first preceding
dot (if the ray precedes all dots we set the value of the function here to $-1$).
We point out that $D-R$ coincides with $\wt  D-\wt R$ on
internal dots and on empty rays\footnote{This is straightforward apart from the case of
an internal dot
whose angular coordinate is not in $(2\pi/d)\NN$, where it holds because, there,
$\wt R$ has been defined as the value of the next ray.}.

The claim follows. An empty ray necessarily has angular coordinate $(2\pi)l/d$ and corresponds
to the sector $\Hcal_{g {J}^{-l}}$.
Since $-\sum_j \lfloor s_j-tw\rfloor$ equals $\wt D+1$,
the degree $(1/2)\deg_{\RW}$
equals
\begin{multline*}
a(g{J}^{-l})-1= \sum_{j=1}^N s_j +\left(- \sum_{j=1}^N \left\lfloor s_j - \frac ld w_j
\right\rfloor\right) -\left(-\left\lfloor -\frac ld d \right\rfloor\right)-1=\\
\sum_j s_j+(\wt D +1)- \wt R-1= \sum_j s_j+ D-  R.
\end{multline*}
On the other hand, for internal dots, the only interesting check concerns the
first dot of one ray $\{\rho \exp(2\pi \cxi t)\mid \rho \in \RR^+\}$. There, the
identities
$\wt R=-\lfloor -td\rfloor$ and $\wt D=-\sum_j \lfloor s_j-tw_j\rfloor$ hold. Therefore
the degree $(1/2)\deg_{\CR}$ of the fundamental class of $\{W_{g \bar \la}=0\}/G\CC^\times$ for
$\la=\exp(-t)$ equals
\begin{multline*}
a(g\bar \la)= \sum_{j=1}^N s_j +\left(- \sum_{j=1}^N \lfloor s_j - t w_j
\rfloor\right) -\left(-\lfloor -t d \rfloor\right)=\\
\sum_j s_j+\wt D- \wt R= \sum_j s_j+ D-  R.
\end{multline*}
\end{proof}

\medskip
\paragraph{Step 5: \emph{The correspondence}.} We finally establish the bidegree
preserving isomorphism. We will be guided by the above diagram which highlights
sets of generators of $H_{\RW}$ (the rays) and sets of generators of $H_{\CR}$ (the dots).
They correspond to each other in a degree-preserving way.

Let us first remark
that the subspaces corresponding to extremal dots in the $\CR$-cohomology
are isomorphic to the subspaces corresponding to the non-empty rays in the
$\RW$-state space. First, if the  angular coordinate of the ray is not contained
in $(2\pi/d)\NN$, then no sector of $H_{\RW}$ is attached to this ray. On the other hand
the primitive cohomology
corresponding to the extremal point on this ray is $\{0\}$ because
the sector is the quotient of a weighted projective stack by a finite group action,
see \eqref{eq:dichotomy}.
Let us focus on a ray $\{\rho \nu \mid \rho \in \RR^+\}$ with $\nu\in \pmmu_d$.
In this case, the extremal dot is the primitive cohomology of the quotient of
a hypersurface
inside a weighted projective stack; this has already been expressed in terms of
$G$-invariant cohomology classes in relative cohomology. Remark \ref{rem:Hdgfilt} yields
the required bidegree-preserving isomorphism.

We finally need to match the internal dots with the empty rays. As remarked above, these
objects correspond to $(p,p)$-classes in the respective $H_{\CR}$ and $H_{\RW}$ spaces
(hyperplane sections and Neveu--Schwarz sectors).
By Lemma \ref{lem:normal}, we only need to provide an involution exchanging
internal dots and empty rays and preserving $D-R$. This is constructed in the next lemma.
\begin{lem}
There exists a $1$-to-$1$ correspondence between internal dots and empty rays
that preserves
$$F=D-R.$$
\end{lem}
\begin{proof}The domain formed by all rays and dots introduced in Step 3 is
totally ordered. The last element is a dot and the first is the real-axis ray $\RR^+$.
Using this order, for any element $n$ different from the last dot
$n+1$ will denote the next element, whereas
for any element $n$  different from the real-axis ray $\RR^+$
we will write $n-1$ for the preceding element.

On the one hand, $n$ is a ray if and only if $F(n-1)=F(n) +1$ or $n=\RR^+$.
On the other hand, $n$ is a dot if and only if $F(n-1)=F(n) -1$.
In other words $F$ is decreasing when it reaches a marking and
is increasing when it reaches a ray. It never varies by more than
$1$. Furthermore the CY condition ensures that $F$ vanishes on
the last value of its domain (in other words the number of
dots equals the number of rays).
It follows that
$F$ may be regarded as a function defined on a set of elements forming a circuit
where the last dot is followed by the first real-axis ray $\RR^+$.
Now notice that if $F$
attains a given value at a given number of internal markings
(going down) it must attains the same value at the same number of
empty rays (going up). Notice that extremal dots and nonempty
rays are the relative maxima and minima of $F$, respectively.
\end{proof}
This completes the proof of Theorem \ref{thm:main}. \qed

\noindent\emph{Proof of Theorem \ref{cor:MS}.}
By Krawitz's main theorem we have $h^{p,q}(W,G;\CC)=h^{N-2-p,q}(W^\trans,G^\trans;\CC)$ (see
\cite[\S2.4]{Kr} and use the fact that $\hat{c}=N-2$). In this way
Theorem \ref{thm:main} yields the claim.
\qed

\begin{rem}\label{rem:norbehaviour}
For $G\subseteq SL_W$, the action of $G\CC^\times$ on $\{W=0\}\subset \CC^N$
satisfies the following property. Consider the point $\pmb x$ in $\{W=0\}$ and
any element $\ga=g\bar \la$ of $G\CC^\times$ fixing $\pmb x$; then, the $(N-1)$-dimensional
representation  $\ga$ in $GL(T_{\pmb x}\{W=0\})$ has determinant $1$.
This happens because $\ga$ acts on the line through $\pmb x$ orthogonal to
$T_{\pmb x}\{W=0\}$ as $z\mapsto \la^d z$. Therefore we have
$\det (\ga\in GL(T_{\pmb x}\{W=0\}) \la^d=\prod_{j=1}^N (g_j\la^{w_j})$; by
the CY condition and $G\subseteq SL_W$, we obtain
$$\det (\ga\in GL(T_{\pmb x}\{W=0\})=\la^{-d} \prod_{j=1}^N (g_j\la^{w_j})=\la^{\sum_j{w_j}-d}\prod_j g_j=1.$$

As a consequence the quotient stack $[X_W/\wt G]$ has no nontrivial orbifold behaviour in codimension
$1$. Therefore, we can relate the ordinary cohomology of the coarse space to
the Chen--Ruan orbifold cohomology of the stack.
Let us assume that the coarse space of $[X_W/\wt G]$, the scheme-theoretic
quotient $X_W/\wt G$, admits a crepant resolution $Z$.
Then, there is a bidegree preserving isomorphism between the cohomology
of $Z$ and the orbifold Chen--Ruan cohomology of $X_W/\wt G$.
In this way Corollaries \ref{cor:coarsecohom} and \ref{cor:LG/CYscheme} follow.
\end{rem}

\section{Examples}\label{sect:exa}

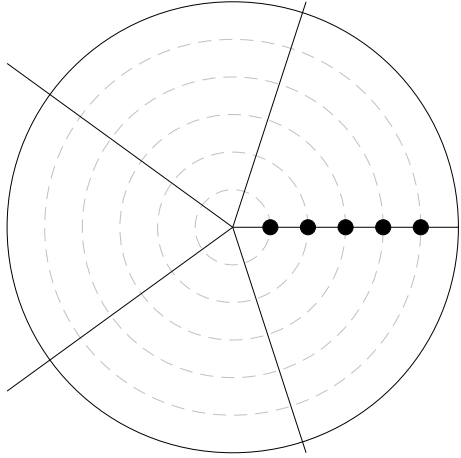
\begin{figure}[t]
% TeXgraph version 1.94 beta-7.5
\begin{tikzpicture}%
\useasboundingbox (-3.,-3.)--(3.,3.);
\pgfsetroundjoin%
%objet1  (Cercle)
\pgfsetstrokecolor{rgb,1:red,0.7529;green,0.7529;blue,0.7529}
\pgfsetlinewidth{0.2pt}
\pgfsetdash{{5pt}{3pt}}{0pt}
\pgfellipse[stroke]{\pgfxy(0,0)}{\pgfxy(1,0)}{\pgfxy(0,1)}
%objet2  (Cercle)
\pgfellipse[stroke]{\pgfxy(0,0)}{\pgfxy(.5,0)}{\pgfxy(0,.5)}
%objet3  (Cercle)
\pgfellipse[stroke]{\pgfxy(0,0)}{\pgfxy(1.5,0)}{\pgfxy(0,1.5)}
%objet4  (Cercle)
\pgfellipse[stroke]{\pgfxy(0,0)}{\pgfxy(2,0)}{\pgfxy(0,2)}
\pgfellipse[stroke]{\pgfxy(0,0)}{\pgfxy(2.5,0)}{\pgfxy(0,2.5)}
\pgfellipse{\pgfxy(0,0)}{\pgfxy(3,0)}{\pgfxy(0,3)}
%\pgfcircle[fillstroke]{\pgfxy(0.5,0)}{.1cm}
%\pgfcircle[fillstroke]{\pgfxy(1.5,0)}{.1cm}
%\pgfcircle[fillstroke]{\pgfxy(1,0)}{.1cm}
%\pgfcircle[fillstroke]{\pgfxy(2.5,0)}{.1cm}
%\pgfcircle[fillstroke]{\pgfxy(2,0)}{.1cm}%Ddroite5  (Utilisateur)
\pgfsetstrokecolor{black}
\pgfsetdash{}{0pt}
\pgfxyline(0,0)(3,0)
\pgfxyline(0,0)(.9747590868, 3)
\pgfxyline(0,0)(-3, 2.179627584  )
\pgfxyline(0,0)(-3 ,-2.179627584  )
\pgfxyline(0,0)(.9747590868 , -3  )
%Ddroite6  (Utilisateur)
\pgfcircle[fillstroke]{\pgfxy(1,0)}{.1cm}
\pgfcircle[fillstroke]{\pgfxy(2,0)}{.1cm}
\pgfcircle[fillstroke]{\pgfxy(.5,0)}{.1cm}
\pgfcircle[fillstroke]{\pgfxy(1.5,0)}{.1cm}
\pgfcircle[fillstroke]{\pgfxy(2.5,0)}{.1cm}
\end{tikzpicture}%
\caption{Diagram of the Fermat quintic in $\PP^4$.\label{fig:quintic}}
\end{figure}
We now recover the examples treated in Section \ref{sect:warmup} and see how they fit in
the diagram illustrated in the course of the proof.
\begin{exa}
Let us consider the case of a degree-$d$ hypersurface in $\PP^{d-2}$ (Example \ref{exa:homog}).
In general, the diagram has $d-1$ empty rays and $d-1$
dots on the real-axis ray.
The diagram for the quintic polynomial in five variables looks as in Figure \ref{fig:quintic}.
The four internal points are the hyperplane sections of the quintic
hypersurface whereas the four empty rays are the Neveu--Schwarz sectors of
the $\RW$ state space. They correspond to each other and the degrees match (they can
be computed following the definition or evaluating the function $D-R$ as in Lemma \ref{lem:normal}
using the diagram).\end{exa}

\begin{exa}\label{exa:K3surf}
Here we illustrate the model in the case of a K3 surface inside a
Gorenstein weighted projective stack. We take the same polynomial as in Example \ref{exa:1inG}, and we get
the diagram
found by Boissi\`ere, Mann and Perroni without
modifications.
In fact, in \cite{BMP},
this diagram is used to describe the sectors of the weighted projective stack
$\PP(2,3,3,4)$; indeed, the dotted rays correspond to the sectors, and the number of
dots lying on one ray corresponds to the dimension of the cohomology of the
corresponding  sector (which,
in turn, is a weighted projective stack).
If we consider the hypersurface where
$W(x_1,\dots,x_4)=x_1^6+x_2^4+x_3^4+x_4^3$
vanishes we can use the diagram as described in Step 3 of the proof.
\begin{figure}[t]
% TeXgraph version 1.94 beta-7.5
\begin{tikzpicture}%
\useasboundingbox (-5.,-5.)--(5.,5.);
\pgfsetroundjoin%
%objet1  (Cercle)
\pgfsetstrokecolor{rgb,1:red,0.7529;green,0.7529;blue,0.7529}
\pgfsetlinewidth{0.2pt}
\pgfsetdash{{5pt}{3pt}}{0pt}
\pgfellipse[stroke]{\pgfxy(0,0)}{\pgfxy(1,0)}{\pgfxy(0,1)}
%objet2  (Cercle)
\pgfellipse[stroke]{\pgfxy(0,0)}{\pgfxy(2,0)}{\pgfxy(0,2)}
%objet3  (Cercle)
\pgfellipse[stroke]{\pgfxy(0,0)}{\pgfxy(3,0)}{\pgfxy(0,3)}
%objet4  (Cercle)
\pgfellipse[stroke]{\pgfxy(0,0)}{\pgfxy(4,0)}{\pgfxy(0,4)}
\pgfellipse{\pgfxy(0,0)}{\pgfxy(5,0)}{\pgfxy(0,5)}
%Ddroite5  (Utilisateur)
\pgfsetstrokecolor{black}
\pgfsetdash{}{0pt}
\pgfxyline(0,0)(5,0)
%Ddroite6  (Utilisateur)
\pgfxyline(0,0)(-5,0)
%Ddroite7  (Utilisateur)
\pgfxyline(0,0)(0,5)
%Ddroite8  (Utilisateur)
\pgfxyline(0,0)(0,-5)
%Ddroite9  (Utilisateur)
\pgfxyline(0,0)(-2.8868,-5)
\pgfxyline(0,0)(2.8868,5)
\pgfxyline(0,0)(5,2.8868)
\pgfxyline(0,0)(-5,-2.8868)
\pgfxyline(0,0)(5,-2.8868)
\pgfxyline(0,0)(-5,2.8868)
\pgfxyline(0,0)(2.8868,-5)
%Ddroite10  (Utilisateur)
\pgfxyline(0,0)(-2.8868,5)
%LabA  (Utilisateur)
\pgfsetfillcolor{black}
\pgfcircle[fillstroke]{\pgfxy(1,0)}{.1cm}
\pgfputat{\pgfxy(0.8232,-0.09)}{\pgftext[right,top]{\color{black}\small 0}}\pgfstroke
%LabB  (Utilisateur)
\pgfcircle[fillstroke]{\pgfxy(2,0)}{.1cm}
\pgfputat{\pgfxy(1.8232,-0.09)}{\pgftext[right,top]{\color{black}\small 1}}\pgfstroke
%LabC  (Utilisateur)
\pgfcircle[fillstroke]{\pgfxy(3,0)}{.1cm}
\pgfputat{\pgfxy(2.8232,-0.09)}{\pgftext[right,top]{\color{black}\small 2}}\pgfstroke
%LabD  (Utilisateur)
\pgfcircle[fillstroke]{\pgfxy(4,0)}{.1cm}
\pgfputat{\pgfxy(3.8232,-0.09)}{\pgftext[right,top]{\color{black}\small 3}}\pgfstroke
%LabE  (Utilisateur)
\pgfcircle[fillstroke]{\pgfxy(0,4)}{.1cm}
\pgfputat{\pgfxy(-0.1768,3.8232)}{\pgftext[right,top]{\color{black}\small 1}}\pgfstroke
%LabF  (Utilisateur)
\pgfcircle[fillstroke]{\pgfxy(-1,1.7321)}{.1cm}
\pgfputat{\pgfxy(-1.1768,1.5553)}{\pgftext[right,top]{\color{black}\small 1}}\pgfstroke
%LabJ  (Utilisateur)
\pgfcircle[fillstroke]{\pgfxy(-1,-1.7321)}{.1cm}
\pgfputat{\pgfxy(-0.8232,-1.9088)}{\pgftext[left,top]{\color{black}\small 1}}\pgfstroke
%LabH  (Utilisateur)
\pgfcircle[fillstroke]{\pgfxy(-1,-0)}{.1cm}
\pgfputat{\pgfxy(-1.1768,-0.1768)}{\pgftext[right,top]{\color{black}\small 1}}\pgfstroke
%LabK  (Utilisateur)
\pgfcircle[fillstroke]{\pgfxy(-1.5,-2.5981)}{.1cm}
\pgfputat{\pgfxy(-1.3232,-2.7749)}{\pgftext[left,top]{\color{black}\small 2}}\pgfstroke
%LabI  (Utilisateur)
\pgfcircle[fillstroke]{\pgfxy(-4,-0)}{.1cm}
\pgfputat{\pgfxy(-4.1768,-0.1768)}{\pgftext[right,top]{\color{black}\small 2}}\pgfstroke
%LabL  (Utilisateur)
\pgfcircle[fillstroke]{\pgfxy(-0,-4)}{.1cm}
\pgfputat{\pgfxy(-0.1768,-4.1768)}{\pgftext[right,top]{\color{black}\small 2}}\pgfstroke
%LabG  (Utilisateur)
\pgfcircle[fillstroke]{\pgfxy(-1.5,2.5981)}{.1cm}
\pgfputat{\pgfxy(-1.6768,2.4213)}{\pgftext[right,top]{\color{black}\small 2}}\pgfstroke
\end{tikzpicture}%
\caption{Diagram of $\{x_1^6+x_2^4+x_3^4+x_4^3=0\}$ inside $\PP(2,3,3,4)$.\label{fig:K3inG}}
\end{figure}
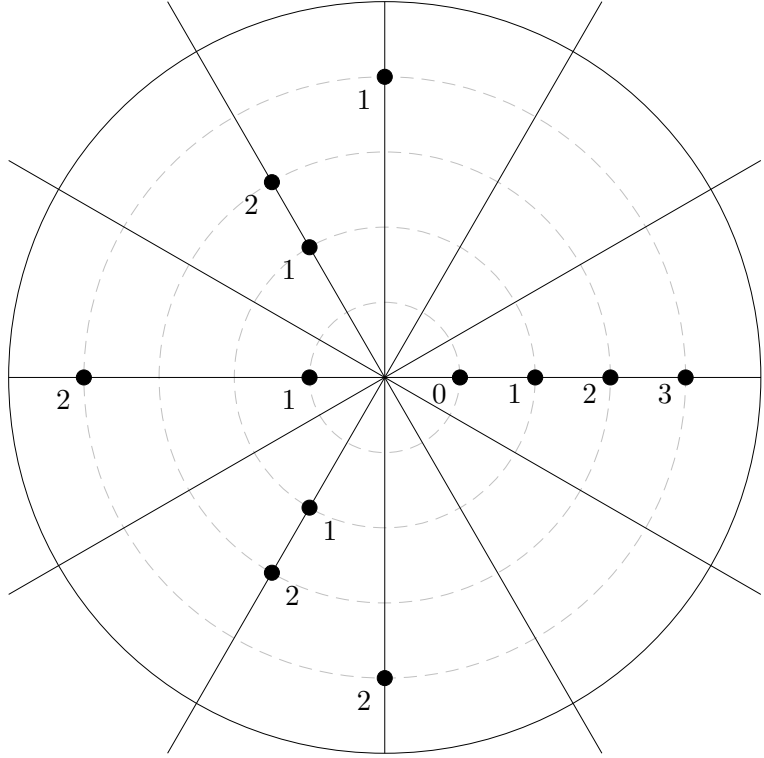
The sectors should be regarded as hypersurfaces lying inside the sectors of the
ambient weighted projective stack. In the surface above
we actually have six dotted rays corresponding to the sectors of the ambient projective stack.
When the ray carries a single dot, the hypersurface is empty. When the ray carries two dots
the hypersurface is $0$-dimensional. Hence, in the example there are only four nonempty  sectors
corresponding to ${J}^0=1, {J}^{-4}, {J}^{-6},$ and ${J}^{-8}$. In general $n$ dots on one ray correspond to an
$(n-2)$-dimensional hypersurface: the first $n-1$ dots counting from the origin are the classes
cut out by $\mathsf 1, \mathsf h, \dots, \mathsf h^{n-2}$, whereas the extremal
dot corresponds to the contribution from primitive cohomology.
Beside each dot we mark the value of $D-R$;
the
reader may check that this coincides with half $\deg_{\CR}$ of the corresponding
class in Chen--Ruan orbifold cohomology (see Example \ref{exa:1inG}).

We leave to the reader
the three-fold
$x_1^{20}+x_2^6+x_3^5+x_4^4+x_5^3$ inside the
Gorenstein weighted projective stacks; we only provide
the combinatorial diagram (see Figure \ref{fig:bigCY} at the end).
%Again, the LG/CY correspondence is based on the fact that the internal dots are as many
%as the empty rays and correspond to each other (CY condition).
\end{exa}

\begin{exa}\label{exa:diagK3nonG}
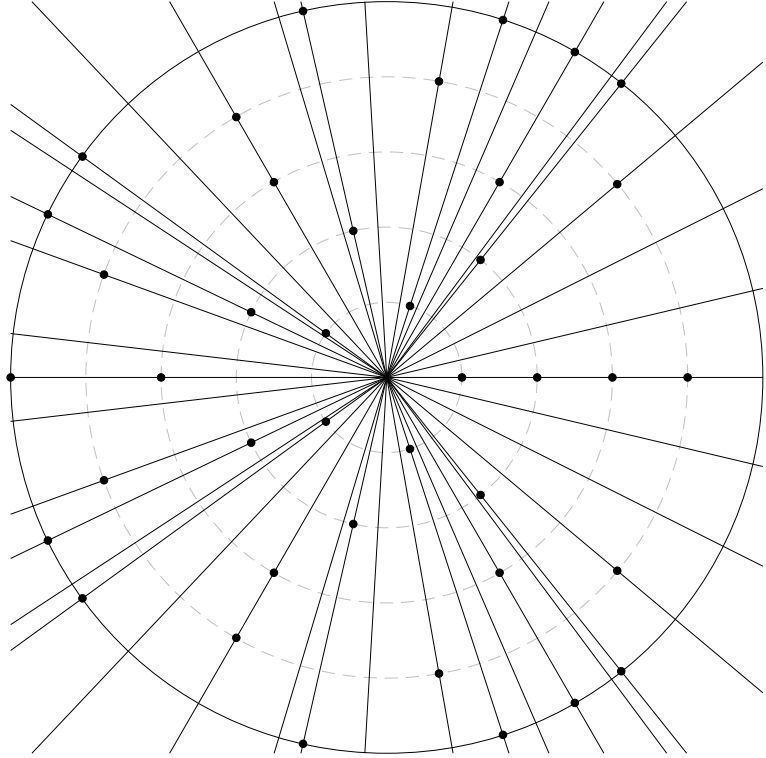
\begin{figure}[t]
% TeXgraph version 1.94 beta-7.5
\begin{tikzpicture}%
\useasboundingbox (-5.,-5.)--(5.,5.);
\pgfsetroundjoin%
%objet1  (Cercle)
\pgfsetstrokecolor{rgb,1:red,0.7529;green,0.7529;blue,0.7529}
\pgfsetlinewidth{0.2pt}
\pgfsetdash{{5pt}{3pt}}{0pt}
\pgfellipse[stroke]{\pgfxy(0,0)}{\pgfxy(1,0)}{\pgfxy(0,1)}
%objet2  (Cercle)
\pgfellipse[stroke]{\pgfxy(0,0)}{\pgfxy(2,0)}{\pgfxy(0,2)}
%objet3  (Cercle)
\pgfellipse[stroke]{\pgfxy(0,0)}{\pgfxy(3,0)}{\pgfxy(0,3)}
%objet4  (Cercle)
\pgfellipse[stroke]{\pgfxy(0,0)}{\pgfxy(4,0)}{\pgfxy(0,4)}
\pgfellipse{\pgfxy(0,0)}{\pgfxy(5,0)}{\pgfxy(0,5)}
%Ddroite5  (Utilisateur)
\pgfsetstrokecolor{black}
\pgfsetdash{}{0pt}
\pgfxyline(0,0)(5,0)
\pgfxyline(0,0)(5,1.185021769)
\pgfxyline(0,0)(5, 2.511094380)
\pgfxyline(0,0)(5,4.195498156)
\pgfxyline(0,0)(3.722362082,5)
\pgfxyline(0,0)(2.156789468,5)
\pgfxyline(0,0)(0.8816349015,5)
\pgfxyline(0,0)(3.987366944,5)
\pgfxyline(0,0)(2.886751345,5)
\pgfxyline(0,0)(1.624598478,5)%here I would stop rays
\pgfxyline(0,0)(-0.2912168377,5)
\pgfxyline(0,0)(-1.496901736,5)
\pgfxyline(0,0)(-2.886751347,5)
\pgfxyline(0,0)(-4.717256704,5)
\pgfxyline(0,0)(-5,3.288551734)
\pgfxyline(0,0)(-5,1.819851172)
                   \pgfxyline(0,0)(-5,0.5844161840)
\pgfxyline(0,0)(-5,3.632712640)
\pgfxyline(0,0)(-1.141217370,5)
\pgfxyline(0,0)(-5, 2.407873095)
\pgfxyline(0,0)(5,-1.185021769)
\pgfxyline(0,0)(5,- 2.511094380)
\pgfxyline(0,0)(5,-4.195498156)
\pgfxyline(0,0)(3.722362082,-5)
\pgfxyline(0,0)(2.156789468,-5)
\pgfxyline(0,0)(0.8816349015,-5)
\pgfxyline(0,0)(3.987366944,-5)
\pgfxyline(0,0)(2.886751345,-5)
\pgfxyline(0,0)(1.624598478,-5)
\pgfxyline(0,0)(-0.2912168377,-5)
\pgfxyline(0,0)(-1.496901736,-5)
\pgfxyline(0,0)(-2.886751347,-5)
\pgfxyline(0,0)(-4.717256704,-5)
\pgfxyline(0,0)(-5,-3.288551734)
\pgfxyline(0,0)(-5,-1.819851172)
                   \pgfxyline(0,0)(-5,-0.5844161840)
\pgfxyline(0,0)(-5,-3.632712640)
\pgfxyline(0,0)(-1.141217370,-5)
\pgfxyline(0,0)(-5, -2.407873095)
\pgfxyline(0,0)(-5,0)
\pgfsetfillcolor{black}
\pgfcircle[fillstroke]{\pgfxy(1,0)}{.05cm}
\pgfcircle[fillstroke]{\pgfxy(2,0)}{.05cm}
\pgfcircle[fillstroke]{\pgfxy(3,0)}{.05cm}
\pgfcircle[fillstroke]{\pgfxy(4,0)}{.05cm}
\pgfcircle[fillstroke]{\pgfxy(1.246979604 , 1.563662965 )}{.05cm}
\pgfcircle[fillstroke]{\pgfxy(-0.4450418678 , 1.949855824 )}{.05cm}
\pgfcircle[fillstroke]{\pgfxy( -1.801937736 , 0.8677674784 )}{.05cm}
\pgfcircle[fillstroke]{\pgfxy(1.246979604 , -1.563662965 )}{.05cm}
\pgfcircle[fillstroke]{\pgfxy(-0.4450418678 , -1.949855824 )}{.05cm}
\pgfcircle[fillstroke]{\pgfxy( -1.801937736 , -0.8677674784 )}{.05cm}
\pgfcircle[fillstroke]{\pgfxy(0.3090169938 , 0.9510565165 )}{.05cm}
\pgfcircle[fillstroke]{\pgfxy(-0.8090169944, 0.5877852522)}{.05cm}
\pgfcircle[fillstroke]{\pgfxy(0.3090169938 , -0.9510565165 )}{.05cm}
\pgfcircle[fillstroke]{\pgfxy(-0.8090169944, -0.5877852522)}{.05cm}
\pgfcircle[fillstroke]{\pgfxy(1.5,2.598076212)}{.05cm}
\pgfcircle[fillstroke]{\pgfxy(-1.5, 2.598076212 )}{.05cm}
\pgfcircle[fillstroke]{\pgfxy(-3,0)}{.05cm}
\pgfcircle[fillstroke]{\pgfxy(2.500000000,4.330127020)}{.05cm}
\pgfcircle[fillstroke]{\pgfxy(-5,0)}{.05cm}
\pgfcircle[fillstroke]{\pgfxy(2.500000000,-4.330127020)}{.05cm}
\pgfcircle[fillstroke]{\pgfxy(1.5, -2.598076212)}{.05cm}
\pgfcircle[fillstroke]{\pgfxy(-1.5, -2.598076212 )}{.05cm}
\pgfcircle[fillstroke]{\pgfxy(3.064177772, 2.571150439)}{.05cm}
\pgfcircle[fillstroke]{\pgfxy(0.6945927092, 3.939231012)}{.05cm}
\pgfcircle[fillstroke]{\pgfxy(-2.000000000, 3.464101616)}{.05cm}
\pgfcircle[fillstroke]{\pgfxy(-3.758770484, 1.368080570)}{.05cm}
\pgfcircle[fillstroke]{\pgfxy(3.064177772, -2.571150439)}{.05cm}
\pgfcircle[fillstroke]{\pgfxy(0.6945927092, -3.939231012)}{.05cm}
\pgfcircle[fillstroke]{\pgfxy(-2.000000000, -3.464101616)}{.05cm}
\pgfcircle[fillstroke]{\pgfxy(-3.758770484, -1.368080570)}{.05cm}
\pgfcircle[fillstroke]{\pgfxy(3.117449009 , 3.909157412)}{.05cm}
\pgfcircle[fillstroke]{\pgfxy(-1.112604670, 4.874639561)}{.05cm}
\pgfcircle[fillstroke]{\pgfxy(-4.504844340, 2.169418696)}{.05cm}
\pgfcircle[fillstroke]{\pgfxy(3.117449009 , -3.909157412)}{.05cm}
\pgfcircle[fillstroke]{\pgfxy(-1.112604670, -4.874639561)}{.05cm}
\pgfcircle[fillstroke]{\pgfxy(-4.504844340, -2.169418696)}{.05cm}
\pgfcircle[fillstroke]{\pgfxy(1.545084969 , 4.755282582)}{.05cm}
\pgfcircle[fillstroke]{\pgfxy(-4.045084972,2.938926261 )}{.05cm}
\pgfcircle[fillstroke]{\pgfxy(1.545084969 , -4.755282582)}{.05cm}
\pgfcircle[fillstroke]{\pgfxy(-4.045084972, -2.938926261 )}{.05cm}
\end{tikzpicture}%
\caption{Diagram of $\{x_1^4x_2+x_2^3x_3+x_3^3x_4+x_4^3=0\}$ inside $\PP(5,7,6,9)$: $40$ rays and $40$ dots.\label{fig:K3nonG}}
\end{figure}
We now illustrate by means of the diagram the case where the
hypersurface is embedded in a nonGorenstein weighted projective stack.
Consider the K3 surface of Example \ref{exa:1nonG}. We illustrate the corresponding
diagram (Figure \ref{fig:K3nonG}).

Two groups should be considered. On the one hand the union of the
roots of unity of order $5,7,6,$ and $9$ (the weights): $H_1=\pmmu_5\cup \pmmu_7 \cup \pmmu_6\cup \pmmu 9$.
On the other hand the roots of unity of order $d=27$ (the degree): $H_2=\pmmu_{27}$. The nonGorenstein
case is characterized by the following feature: $H_2\not \subseteq H_1$.

Let us now go through the definition. We draw a ray for every element of $H_1\cup H_2$. In this way
we have $40$ rays ($13$ of them are special because they correspond to elements of
$H_2\setminus H_1$. We mark dots on the four circles corresponding to the four coordinates:
$5$ dots on the first, $7$ dots on the second, $6$ on the third, and $9$ on the fourth.
Following the construction of Step 3 of the proof, we mark 13 further dots with polar coordinate $N+1$.

The presence of rays whose angular coordinate is not in
$2\pi\cxi\{0,\frac1{27},\dots,\frac{26}{27}\}$
corresponds to the fact that there are sectors that do not intersect transversely $\{W=0\}$.
The correspondence still holds because the presence of extra rays is balanced by the presence of extra dots.
\end{exa}

\begin{exa}\label{exa:demo}
This example is meant to illustrate the setup of the proof in the more delicate cases
where nontrivial $\langle {J}\rangle$-cosets are involved.
We consider the cubic equation already studied in Example \ref{exa:groupquot}, \emph{i.e.}
$x_1^2x_2+x_2^2x_3+x_3^3=0$, and the order-12 cyclic group $G=\Aut(W)$.

As in the proof, we proceed coset by coset.
Note that $\ga^4={J}$, therefore the
natural choices corresponding to
$g^{(1)}, g^{(2)}, g^{(3)}, g^{(4)}$ in the proof are
$\ga^0,\ga^1, \ga^2,\ga^3$.

We start from the coset attached to $g=\ga^0=(1,1,1)$ and we apply the previous construction.
The terms $(g_1, \dots, g_N)$
are the $N$ coordinates of $g\in (\CC^\times)^N$: in
this case they are all equal to $1$.
We have $\{\al\mid \al^{w_j}=g_j\}=\{1\}$ because the weights are all equal to $1$.
We have
$$\pmmu_d\cup \bigcup_{j=1}^N \left\{\al\mid \al^{w_j}=g_j\right\}=\pmmu_d,$$
hence there are three rays (as many as $d$, which equals $3$).
Similarly there are three dots, as many as the solutions (in the variables
$\nu$ and $j$) of $\nu^{w_j}=1$: $(\nu,j)$ is necessarily $(1,1), (1,2)$, or $(1,3)$.
Note that the further dots mentioned in the construction of
the model do not occur in this coset because
$\bigcup_{j=1}^N \{\al\mid \al^{w_j}=g_j\}$ is contained in $\pmmu_d$.
The picture is that of Figure \ref{fig:1stcoset}.
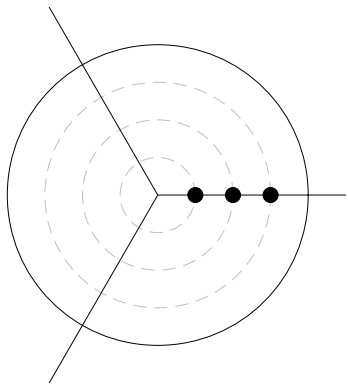
\begin{figure}[h]
% TeXgraph version 1.94 beta-7.5
\begin{tikzpicture}%
\useasboundingbox (-2.5,-2.5)--(2.5,2.5);
\pgfsetroundjoin%
%objet1  (Cercle)
\pgfsetstrokecolor{rgb,1:red,0.7529;green,0.7529;blue,0.7529}
\pgfsetlinewidth{0.2pt}
\pgfsetdash{{5pt}{3pt}}{0pt}
\pgfellipse[stroke]{\pgfxy(0,0)}{\pgfxy(.5,0)}{\pgfxy(0,.5)}
%objet2  (Cercle)
\pgfellipse[stroke]{\pgfxy(0,0)}{\pgfxy(1,0)}{\pgfxy(0,1)}
%objet3  (Cercle)
\pgfellipse[stroke]{\pgfxy(0,0)}{\pgfxy(1.5,0)}{\pgfxy(0,1.5)}
%objet4  (Cercle)
\pgfellipse{\pgfxy(0,0)}{\pgfxy(2,0)}{\pgfxy(0,2)}
%Ddroite5  (Utilisateur)
\pgfsetstrokecolor{black}
\pgfsetdash{}{0pt}
\pgfxyline(0,0)(2.5,0)
%Ddroite6  (Utilisateur)
%\pgfxyline(0,0)(-5,0)
%Ddroite7  (Utilisateur)
%\pgfxyline(0,0)(0,5)
%Ddroite8  (Utilisateur)
%\pgfxyline(0,0)(0,-5)
%Ddroite9  (Utilisateur)
\pgfxyline(0,0)(-1.4434,-2.5)
%Ddroite10  (Utilisateur)
\pgfxyline(0,0)(-1.4434,2.5)
%LabA  (Utilisateur)
\pgfsetfillcolor{black}
\pgfcircle[fillstroke]{\pgfxy(.5,0)}{.1cm}
%\pgfputat{\pgfxy(0.4116,-0.1768)}{\pgftext[right,top]{\color{black}\small 0}}\pgfstroke
%LabB  (Utilisateur)
\pgfcircle[fillstroke]{\pgfxy(1,0)}{.1cm}
%\pgfputat{\pgfxy(0.9116,-0.1768)}{\pgftext[right,top]{\color{black}\small 1}}\pgfstroke
%LabC  (Utilisateur)
\pgfcircle[fillstroke]{\pgfxy(1.5,0)}{.1cm}
%\pgfputat{\pgfxy(1.4116,-0.1768)}{\pgftext[right,top]{\color{black}\small 2}}\pgfstroke
%LabD  (Utilisateur)
\end{tikzpicture}%
\caption{Diagram attached to $(1,1,1)$. \label{fig:1stcoset}}
\end{figure}

We can move on to the coset corresponding to $g=\ga$. This time the three coordinates differ
$g_1=\exp(2\pi\cxi 1/12)$: there is a single solution to $\al^{w_1}=g_1$ which is $\al=\exp(2\pi\cxi 1/12)$.
Similarly there is a single solution to $\al^{w_2}=g_2$, \emph{i.e.} $\al$ equal to $\exp(2\pi\cxi 10/12)$,
and
there is a single solution to $\al^{w_3}=g_3$, \emph{i.e.} $\al$ equal to $\exp(2\pi\cxi 4/12)$.
We have
$$\pmmu_d\cup \bigcup_{j=1}^N \left\{\al\mid \al^{w_j}=g_j\right\}=\textstyle{\pmmu_3\cup \left\{\exp(2\pi\cxi \frac 1{12}),
\exp(2\pi\cxi \frac{10}{12})\right\}.}$$
Therefore we draw five rays (whose angular coordinates range among those of the above set).
Following the rules of Section \ref{sect:main} we draw five dots:
$$\textstyle{1\exp(2\pi\cxi \frac1{12}), 2\exp(2\pi\cxi \frac{10}{12}), 3\exp(2\pi\cxi \frac{4}{12}), 4\exp(2\pi\cxi \frac1{12}),
4\exp(2\pi\cxi \frac{10}{12}),}$$
where the last two dots correspond to the set
$\left( \bigcup_{j=1}^N \left\{\al\mid \al^{w_j}=g_j\right\}\right)\setminus \pmmu_d$ which consists of
two elements: $\exp(2\pi\cxi 1/12)$ and $\exp(2\pi\cxi 10/12)$.

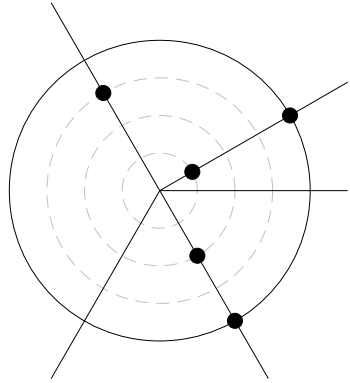
\begin{figure}[h]
% TeXgraph version 1.94 beta-7.5
\begin{tikzpicture}%
\useasboundingbox (-2.5,-2.5)--(2.5,2.5);
\pgfsetroundjoin%
%objet1  (Cercle)
\pgfsetstrokecolor{rgb,1:red,0.7529;green,0.7529;blue,0.7529}
\pgfsetlinewidth{0.2pt}
\pgfsetdash{{5pt}{3pt}}{0pt}
\pgfellipse[stroke]{\pgfxy(0,0)}{\pgfxy(.5,0)}{\pgfxy(0,.5)}
%objet2  (Cercle)
\pgfellipse[stroke]{\pgfxy(0,0)}{\pgfxy(1,0)}{\pgfxy(0,1)}
%objet3  (Cercle)
\pgfellipse[stroke]{\pgfxy(0,0)}{\pgfxy(1.5,0)}{\pgfxy(0,1.5)}
%objet4  (Cercle)
\pgfellipse{\pgfxy(0,0)}{\pgfxy(2,0)}{\pgfxy(0,2)}
%Ddroite5  (Utilisateur)
\pgfsetstrokecolor{black}
\pgfsetdash{}{0pt}
\pgfxyline(0,0)(2.5,0)
%Ddroite6  (Utilisateur)
%\pgfxyline(0,0)(-5,0)
%Ddroite7  (Utilisateur)
%\pgfxyline(0,0)(0,5)
%Ddroite8  (Utilisateur)
%\pgfxyline(0,0)(0,-5)
%Ddroite9  (Utilisateur)
\pgfxyline(0,0)(-1.4434,-2.5)
%Ddroite10  (Utilisateur)
\pgfxyline(0,0)(-1.4434,2.5)
%\pgfxyline(0,0)(0,-5)
%Ddroite11  (Utilisateur)
\pgfxyline(0,0)(2.5,1.4434)
%Ddroite12  (Utilisateur)
\pgfxyline(0,0)(1.4434,-2.5)
%LabA  (Utilisateur)
\pgfsetfillcolor{black}
%\pgfcircle[fillstroke]{\pgfxy(.5,0)}{.1cm}
%\pgfputat{\pgfxy(0.4116,-0.1768)}{\pgftext[right,top]{\color{black}\small 0}}\pgfstroke
%LabB  (Utilisateur)
%\pgfcircle[fillstroke]{\pgfxy(1,0)}{.1cm}
%\pgfputat{\pgfxy(0.9116,-0.1768)}{\pgftext[right,top]{\color{black}\small 1}}\pgfstroke
%LabC  (Utilisateur)
%\pgfcircle[fillstroke]{\pgfxy(1.5,0)}{.1cm}
%\pgfputat{\pgfxy(1.4116,-0.1768)}{\pgftext[right,top]{\color{black}\small 2}}\pgfstroke
%LabD  (Utilisateur)
%\pgfcircle[fillstroke]{\pgfxy(2,0)}{.1cm}
%\pgfputat{\pgfxy(1.9116,-0.1768)}{\pgftext[right,top]{\color{black}\small 3}}\pgfstroke
%LabE  (Utilisateur)
\pgfcircle[fillstroke]{\pgfxy(1,-1.732)}{.1cm}
\pgfcircle[fillstroke]{\pgfxy(0.5,-0.866)}{.1cm}
%\pgfputat{\pgfxy(-0.1768,3.8232)}{\pgftext[right,top]{\color{black}\small 4}}\pgfstroke
%LabF  (Utilisateur)
\pgfcircle[fillstroke]{\pgfxy(0.433,0.25)}{.1cm}
\pgfcircle[fillstroke]{\pgfxy(1.732,1)}{.1cm}
\pgfcircle[fillstroke]{\pgfxy(-0.75,1.299)}{.1cm}
%\pgfputat{\pgfxy(-1.1768,1.5553)}{\pgftext[right,top]{\color{black}\small 5}}\pgfstroke
%LabJ  (Utilisateur)
%\pgfcircle[fillstroke]{\pgfxy(-1,-1.7321)}{.1cm}
%\pgfputat{\pgfxy(-0.8232,-1.9088)}{\pgftext[left,top]{\color{black}\small 9}}\pgfstroke
%LabH  (Utilisateur)
%\pgfcircle[fillstroke]{\pgfxy(-1,-0)}{.1cm}
%\pgfputat{\pgfxy(-1.1768,-0.1768)}{\pgftext[right,top]{\color{black}\small 7}}\pgfstroke
%LabK  (Utilisateur)
%\pgfcircle[fillstroke]{\pgfxy(-1.5,-2.5981)}{.1cm}
%\pgfputat{\pgfxy(-1.3232,-2.7749)}{\pgftext[left,top]{\color{black}\small 10}}\pgfstroke
%LabI  (Utilisateur)
%\pgfcircle[fillstroke]{\pgfxy(-4,-0)}{.1cm}
%\pgfputat{\pgfxy(-4.1768,-0.1768)}{\pgftext[right,top]{\color{black}\small 8}}\pgfstroke
%LabL  (Utilisateur)
%\pgfcircle[fillstroke]{\pgfxy(-0,-4)}{.1cm}
%\pgfputat{\pgfxy(-0.1768,-4.1768)}{\pgftext[right,top]{\color{black}\small 11}}\pgfstroke
%LabG  (Utilisateur)
%\pgfcircle[fillstroke]{\pgfxy(-1.5,2.5981)}{.1cm}
%\pgfputat{\pgfxy(-1.6768,2.4213)}{\pgftext[right,top]{\color{black}\small 6}}\pgfstroke
\end{tikzpicture}%
\caption{Diagram for $\exp(2\pi\cxi(1/12,10/12,4/12))$. \label{fig:2ndcoset}}
\end{figure}

The analysis of the third and fourth cosets is completely analogous
to that we just carried out and yields Figures \ref{fig:nastyray} and \ref{fig:4thcoset}.
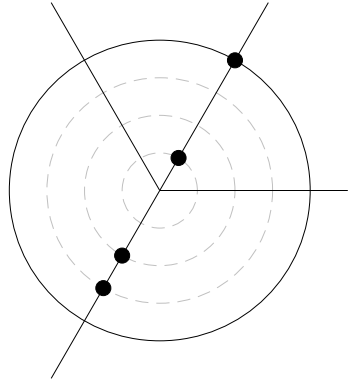
\begin{figure}[h]
% TeXgraph version 1.94 beta-7.5
\begin{tikzpicture}%
\useasboundingbox (-2.5,-2.5)--(2.5,2.5);
\pgfsetroundjoin%
%objet1  (Cercle)
\pgfsetstrokecolor{rgb,1:red,0.7529;green,0.7529;blue,0.7529}
\pgfsetlinewidth{0.2pt}
\pgfsetdash{{5pt}{3pt}}{0pt}
\pgfellipse[stroke]{\pgfxy(0,0)}{\pgfxy(.5,0)}{\pgfxy(0,.5)}
%objet2  (Cercle)
\pgfellipse[stroke]{\pgfxy(0,0)}{\pgfxy(1,0)}{\pgfxy(0,1)}
%objet3  (Cercle)
\pgfellipse[stroke]{\pgfxy(0,0)}{\pgfxy(1.5,0)}{\pgfxy(0,1.5)}
%objet4  (Cercle)
\pgfellipse{\pgfxy(0,0)}{\pgfxy(2,0)}{\pgfxy(0,2)}
%Ddroite5  (Utilisateur)
\pgfsetstrokecolor{black}
\pgfsetdash{}{0pt}
\pgfxyline(0,0)(2.5,0)
%Ddroite6  (Utilisateur)
%\pgfxyline(0,0)(-5,0)
%Ddroite7  (Utilisateur)
%\pgfxyline(0,0)(0,5)
%Ddroite8  (Utilisateur)
%\pgfxyline(0,0)(0,-5)
%Ddroite9  (Utilisateur)
\pgfxyline(0,0)(-1.4434,-2.5)
%Ddroite10  (Utilisateur)
\pgfxyline(0,0)(-1.4434,2.5)
%\pgfxyline(0,0)(0,-5)
%Ddroite11  (Utilisateur)
%\pgfxyline(0,0)(2.5,1.4434)
%Ddroite12  (Utilisateur)
\pgfxyline(0,0)(1.4434,2.5)
%LabA  (Utilisateur)
\pgfsetfillcolor{black}
%\pgfcircle[fillstroke]{\pgfxy(.5,0)}{.1cm}
%\pgfputat{\pgfxy(0.4116,-0.1768)}{\pgftext[right,top]{\color{black}\small 0}}\pgfstroke
%LabB  (Utilisateur)
%\pgfcircle[fillstroke]{\pgfxy(1,0)}{.1cm}
%\pgfputat{\pgfxy(0.9116,-0.1768)}{\pgftext[right,top]{\color{black}\small 1}}\pgfstroke
%LabC  (Utilisateur)
%\pgfcircle[fillstroke]{\pgfxy(1.5,0)}{.1cm}
%\pgfputat{\pgfxy(1.4116,-0.1768)}{\pgftext[right,top]{\color{black}\small 2}}\pgfstroke
%LabD  (Utilisateur)
%\pgfcircle[fillstroke]{\pgfxy(2,0)}{.1cm}
%\pgfputat{\pgfxy(1.9116,-0.1768)}{\pgftext[right,top]{\color{black}\small 3}}\pgfstroke
%LabE  (Utilisateur)
\pgfcircle[fillstroke]{\pgfxy(1,1.732)}{.1cm}
\pgfcircle[fillstroke]{\pgfxy(0.25,0.433)}{.1cm}
\pgfcircle[fillstroke]{\pgfxy(-0.5,-0.866)}{.1cm}
\pgfcircle[fillstroke]{\pgfxy(-0.75,-1.299)}{.1cm}
%\pgfputat{\pgfxy(-0.1768,3.8232)}{\pgftext[right,top]{\color{black}\small 4}}\pgfstroke
%LabF  (Utilisateur)
%\pgfcircle[fillstroke]{\pgfxy(1.732,1)}{.1cm}
%\pgfputat{\pgfxy(-1.1768,1.5553)}{\pgftext[right,top]{\color{black}\small 5}}\pgfstroke
%LabJ  (Utilisateur)
%\pgfcircle[fillstroke]{\pgfxy(-1,-1.7321)}{.1cm}
%\pgfputat{\pgfxy(-0.8232,-1.9088)}{\pgftext[left,top]{\color{black}\small 9}}\pgfstroke
%LabH  (Utilisateur)
%\pgfcircle[fillstroke]{\pgfxy(-1,-0)}{.1cm}
%\pgfputat{\pgfxy(-1.1768,-0.1768)}{\pgftext[right,top]{\color{black}\small 7}}\pgfstroke
%LabK  (Utilisateur)
%\pgfcircle[fillstroke]{\pgfxy(-1.5,-2.5981)}{.1cm}
%\pgfputat{\pgfxy(-1.3232,-2.7749)}{\pgftext[left,top]{\color{black}\small 10}}\pgfstroke
%LabI  (Utilisateur)
%\pgfcircle[fillstroke]{\pgfxy(-4,-0)}{.1cm}
%\pgfputat{\pgfxy(-4.1768,-0.1768)}{\pgftext[right,top]{\color{black}\small 8}}\pgfstroke
%LabL  (Utilisateur)
%\pgfcircle[fillstroke]{\pgfxy(-0,-4)}{.1cm}
%\pgfputat{\pgfxy(-0.1768,-4.1768)}{\pgftext[right,top]{\color{black}\small 11}}\pgfstroke
%LabG  (Utilisateur)
%\pgfcircle[fillstroke]{\pgfxy(-1.5,2.5981)}{.1cm}
%\pgfputat{\pgfxy(-1.6768,2.4213)}{\pgftext[right,top]{\color{black}\small 6}}\pgfstroke
\end{tikzpicture}%
\caption{Diagram for $\exp(2\pi\cxi(2/12,8/12,8/12))$.\label{fig:nastyray}}
\end{figure}

\begin{figure}[h]
% TeXgraph version 1.94 beta-7.5
\begin{tikzpicture}%
\useasboundingbox (-2.5,-2.5)--(2.5,2.5);
\pgfsetroundjoin%
%objet1  (Cercle)
\pgfsetstrokecolor{rgb,1:red,0.7529;green,0.7529;blue,0.7529}
\pgfsetlinewidth{0.2pt}
\pgfsetdash{{5pt}{3pt}}{0pt}
\pgfellipse[stroke]{\pgfxy(0,0)}{\pgfxy(.5,0)}{\pgfxy(0,.5)}
%objet2  (Cercle)
\pgfellipse[stroke]{\pgfxy(0,0)}{\pgfxy(1,0)}{\pgfxy(0,1)}
%objet3  (Cercle)
\pgfellipse[stroke]{\pgfxy(0,0)}{\pgfxy(1.5,0)}{\pgfxy(0,1.5)}
%objet4  (Cercle)
\pgfellipse{\pgfxy(0,0)}{\pgfxy(2,0)}{\pgfxy(0,2)}
%Ddroite5  (Utilisateur)
\pgfsetstrokecolor{black}
\pgfsetdash{}{0pt}
\pgfxyline(0,0)(2.5,0)
%Ddroite6  (Utilisateur)
%\pgfxyline(0,0)(-5,0)
%Ddroite7  (Utilisateur)
%\pgfxyline(0,0)(0,5)
%Ddroite8  (Utilisateur)
%\pgfxyline(0,0)(0,-5)
%Ddroite9  (Utilisateur)
\pgfxyline(0,0)(-1.4434,-2.5)
%Ddroite10  (Utilisateur)
\pgfxyline(0,0)(-1.4434,2.5)
%LabA  (Utilisateur)
\pgfsetfillcolor{black}
\pgfcircle[fillstroke]{\pgfxy(-1,0)}{.1cm}
\pgfxyline(0,0)(-2.5,0)
\pgfcircle[fillstroke]{\pgfxy(-2,0)}{.1cm}
%\pgfputat{\pgfxy(0.4116,-0.1768)}{\pgftext[right,top]{\color{black}\small 0}}\pgfstroke
%LabB  (Utilisateur)
\pgfcircle[fillstroke]{\pgfxy(0,0.5)}{.1cm}
\pgfxyline(0,0)(0,2.5)
\pgfcircle[fillstroke]{\pgfxy(0,2)}{.1cm}
%\pgfputat{\pgfxy(0.9116,-0.1768)}{\pgftext[right,top]{\color{black}\small 1}}\pgfstroke
%LabC  (Utilisateur)
\pgfcircle[fillstroke]{\pgfxy(1.5,0)}{.1cm}
%\pgfputat{\pgfxy(1.4116,-0.1768)}{\pgftext[right,top]{\color{black}\small 2}}\pgfstroke
%LabD  (Utilisateur)
%\pgfcircle[fillstroke]{\pgfxy(2,0)}{.1cm}
%\pgfputat{\pgfxy(1.9116,-0.1768)}{\pgftext[right,top]{\color{black}\small 3}}\pgfstroke
%LabE  (Utilisateur)
%\pgfcircle[fillstroke]{\pgfxy(0,2)}{.1cm}
%\pgfputat{\pgfxy(-0.1768,3.8232)}{\pgftext[right,top]{\color{black}\small 4}}\pgfstroke
%LabF  (Utilisateur)
%\pgfcircle[fillstroke]{\pgfxy(-.5,0.8660)}{.1cm}
%\pgfputat{\pgfxy(-1.1768,1.5553)}{\pgftext[right,top]{\color{black}\small 5}}\pgfstroke
%LabJ  (Utilisateur)
%\pgfcircle[fillstroke]{\pgfxy(-1,-1.7321)}{.1cm}
%\pgfputat{\pgfxy(-0.8232,-1.9088)}{\pgftext[left,top]{\color{black}\small 9}}\pgfstroke
%LabH  (Utilisateur)
%\pgfcircle[fillstroke]{\pgfxy(-1,-0)}{.1cm}
%\pgfputat{\pgfxy(-1.1768,-0.1768)}{\pgftext[right,top]{\color{black}\small 7}}\pgfstroke
%LabK  (Utilisateur)
%\pgfcircle[fillstroke]{\pgfxy(-1.5,-2.5981)}{.1cm}
%\pgfputat{\pgfxy(-1.3232,-2.7749)}{\pgftext[left,top]{\color{black}\small 10}}\pgfstroke
%LabI  (Utilisateur)
%\pgfcircle[fillstroke]{\pgfxy(-4,-0)}{.1cm}
%\pgfputat{\pgfxy(-4.1768,-0.1768)}{\pgftext[right,top]{\color{black}\small 8}}\pgfstroke
%LabL  (Utilisateur)
%\pgfcircle[fillstroke]{\pgfxy(-0,-4)}{.1cm}
%\pgfputat{\pgfxy(-0.1768,-4.1768)}{\pgftext[right,top]{\color{black}\small 11}}\pgfstroke
%LabG  (Utilisateur)
%\pgfcircle[fillstroke]{\pgfxy(-1.5,2.5981)}{.1cm}
%\pgfputat{\pgfxy(-1.6768,2.4213)}{\pgftext[right,top]{\color{black}\small 6}}\pgfstroke
\end{tikzpicture}%
\caption{Diagram for $\exp(2\pi\cxi(3/12,6/12,0))$. \label{fig:4thcoset}}
\end{figure}
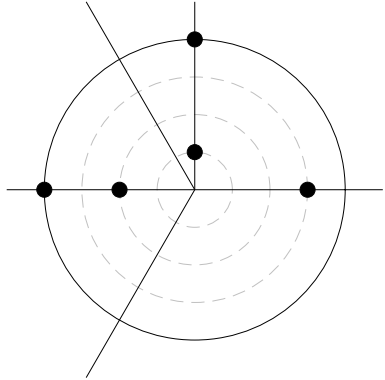

This setting allows one
to check that there is a degree-preserving isomorphism.
We can focus on the eight empty rays (on the $\RW$ side) and
compare them to the eight internal points.
Using Lemma \ref{lem:normal} we get the degrees on the four diagrams.
On Figure \ref{fig:1stcoset} there are two internal dots on the real axis
for which $\deg_{\CR}$ is $0$ and $1$ (if we read in lexicographic order),
and --- correspondingly --- two empty rays for which
$\deg_{\RW}$ is $1$ and $0$ (reading in the sense of the angular coordinate).
It is an
interesting exercise to verify that all the internal dots and empty rays
appearing in Figure \ref{fig:2ndcoset}
have degree $1/2$ (twice $a(h)+D-R$),
all internal dots on Figure \ref{fig:nastyray} have degree $1$, and, finally,
all internal dots on Figure \ref{fig:4thcoset} have degree $3/2$.
This matches the orbifold curve,
Figure \ref{fig:sectorsorbicurve}.

\medskip

{\bf Acknowledgements.\/} Special thanks to Victor Batyrev for helpful conversations on
classical mirror symmetry conjecture in toric geometry.
We also would like to thank Samuel Boissi\`ere, Gavin Brown, Alexandru Dimca,
Bashar Dudin, Arthur Greenspoon,
Kentaro Hori, Marc Krawitz,  Catriona Maclean,
Johannes Nicaise, Jan Nagel, and Matthieu Romagny
for useful discussions.

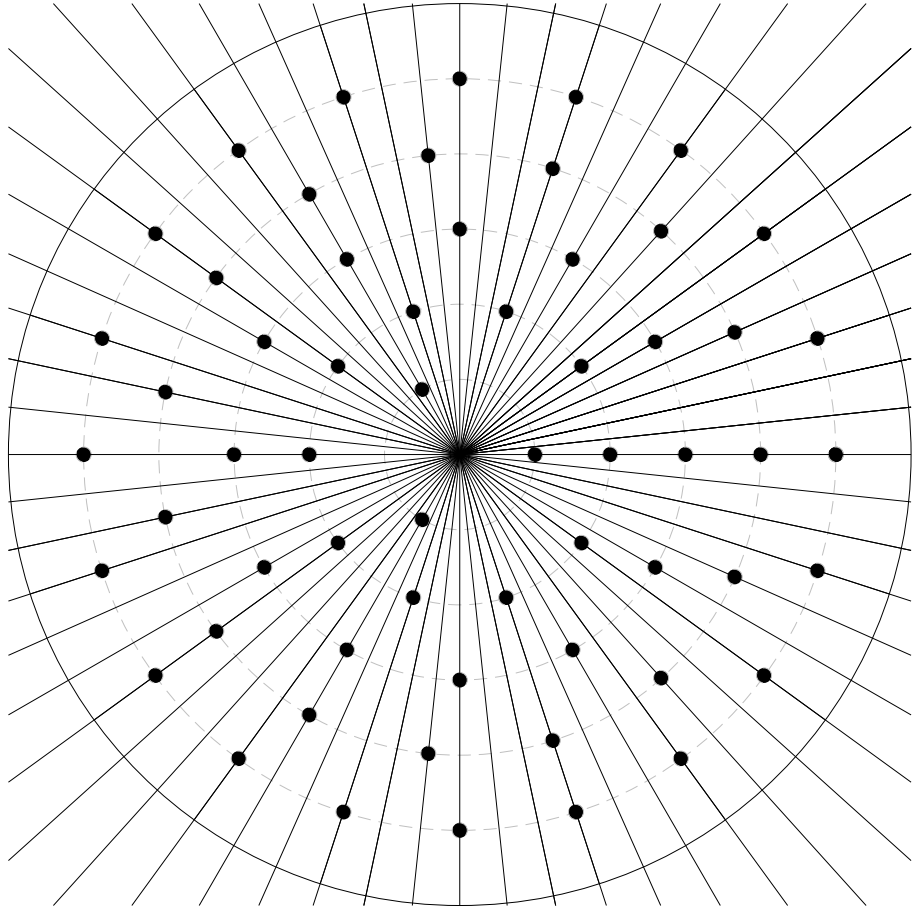
\begin{figure}[h]
% TeXgraph version 1.94 beta-7.5
\begin{tikzpicture}%
\useasboundingbox (-6,-6)--(6,6);
\pgfsetroundjoin%
%objet1  (Cercle)
\pgfsetstrokecolor{rgb,1:red,0.7529;green,0.7529;blue,0.7529}
\pgfsetlinewidth{0.2pt}
\pgfsetdash{{5pt}{3pt}}{0pt}
\pgfellipse[stroke]{\pgfxy(0,0)}{\pgfxy(1,0)}{\pgfxy(0,1)}
\pgfcircle[fillstroke]{\pgfxy(-.5,0.866)}{.1cm}
\pgfcircle[fillstroke]{\pgfxy(-.5,-0.866)}{.1cm}
\pgfcircle[fillstroke]{\pgfxy(1,0)}{.1cm}
%objet2  (Cercle)
\pgfellipse[stroke]{\pgfxy(0,0)}{\pgfxy(2,0)}{\pgfxy(0,2)}
\pgfcircle[fillstroke]{\pgfxy(1.618033989, + 1.175570505)}{.1cm}
\pgfcircle[fillstroke]{\pgfxy(0.6180339876, + 1.902113033)}{.1cm}
\pgfcircle[fillstroke]{\pgfxy(-0.6180339876, + 1.902113033)}{.1cm}
\pgfcircle[fillstroke]{\pgfxy(-1.6180339876, + 1.1755)}{.1cm}
\pgfcircle[fillstroke]{\pgfxy(-2, 0)}{.1cm}
\pgfcircle[fillstroke]{\pgfxy(2,0)}{.1cm}
\pgfcircle[fillstroke]{\pgfxy(-1.618033989, - 1.175570505)}{.1cm}
\pgfcircle[fillstroke]{\pgfxy(-0.6180339876, - 1.902113033)}{.1cm}
\pgfcircle[fillstroke]{\pgfxy(+0.6180339876, - 1.902113033)}{.1cm}
\pgfcircle[fillstroke]{\pgfxy(1.6180339876, - 1.1755)}{.1cm}
%objet3  (Cercle)
\pgfellipse[stroke]{\pgfxy(0,0)}{\pgfxy(3,0)}{\pgfxy(0,3)}
\pgfcircle[fillstroke]{\pgfxy(0, 3)}{.1cm}
\pgfcircle[fillstroke]{\pgfxy(0,-3)}{.1cm}
\pgfcircle[fillstroke]{\pgfxy(2.598, 1.5)}{.1cm}
\pgfcircle[fillstroke]{\pgfxy(1.5, 2.598)}{.1cm}
\pgfcircle[fillstroke]{\pgfxy(3,0)}{.1cm}
\pgfcircle[fillstroke]{\pgfxy(-1.5, 2.598)}{.1cm}
\pgfcircle[fillstroke]{\pgfxy(-2.598, 1.5)}{.1cm}
\pgfcircle[fillstroke]{\pgfxy(-2.598,- 1.5)}{.1cm}
\pgfcircle[fillstroke]{\pgfxy(-1.5, -2.598)}{.1cm}
\pgfcircle[fillstroke]{\pgfxy(-3,0)}{.1cm}
\pgfcircle[fillstroke]{\pgfxy(1.5, -2.598)}{.1cm}
\pgfcircle[fillstroke]{\pgfxy(2.598, -1.5)}{.1cm}
%objet4  (Cercle)
\pgfellipse[stroke]{\pgfxy(0,0)}{\pgfxy(4,0)}{\pgfxy(0,4)}
\pgfcircle[fillstroke]{\pgfxy(4, 0)}{.1cm}
\pgfcircle[fillstroke]{\pgfxy(3.654181831 , 1.626946572 )}{.1cm}
\pgfcircle[fillstroke]{\pgfxy( 2.676522425 ,+ 2.972579302 )}{.1cm}
\pgfcircle[fillstroke]{\pgfxy(1.236067975, + 3.804226066 )}{.1cm}
\pgfcircle[fillstroke]{\pgfxy(-0.4181138536, + 3.978087582 )}{.1cm}
\pgfcircle[fillstroke]{\pgfxy(-2. ,+ 3.464101616 )}{.1cm}
\pgfcircle[fillstroke]{\pgfxy( -3.236067978 ,+ 2.351141009 )}{.1cm}
\pgfcircle[fillstroke]{\pgfxy(-3.912590404 ,+ 0.8316467608  )}{.1cm}
\pgfcircle[fillstroke]{\pgfxy(3.654181831 , -1.626946572 )}{.1cm}
\pgfcircle[fillstroke]{\pgfxy(2.676522425 ,- 2.972579302 )}{.1cm}
\pgfcircle[fillstroke]{\pgfxy(1.236067975, - 3.804226066 )}{.1cm}
\pgfcircle[fillstroke]{\pgfxy(-0.4181138536, - 3.978087582 )}{.1cm}
\pgfcircle[fillstroke]{\pgfxy(-2. ,- 3.464101616 )}{.1cm}
\pgfcircle[fillstroke]{\pgfxy(-3.236067978 ,- 2.351141009 )}{.1cm}
\pgfcircle[fillstroke]{\pgfxy(-3.912590404 ,- 0.8316467608  )}{.1cm}
%objet5  (Cercle)
\pgfellipse[stroke]{\pgfxy(0,0)}{\pgfxy(5,0)}{\pgfxy(0,5)}
\pgfcircle[fillstroke]{\pgfxy(5, 0)}{.1cm}
\pgfcircle[fillstroke]{\pgfxy(4.755282582, + 1.545084972 )}{.1cm}
\pgfcircle[fillstroke]{\pgfxy( 4.045084972, + 2.938926262 )}{.1cm}
\pgfcircle[fillstroke]{\pgfxy(2.938926261 ,+ 4.045084972 )}{.1cm}
\pgfcircle[fillstroke]{\pgfxy(1.545084969 ,+ 4.755282582  )}{.1cm}
\pgfcircle[fillstroke]{\pgfxy(-5. ,0 )}{.1cm}
\pgfcircle[fillstroke]{\pgfxy(0 ,5 )}{.1cm}
\pgfcircle[fillstroke]{\pgfxy(0,-5)}{.1cm}
\pgfcircle[fillstroke]{\pgfxy(-1.545084971 ,+ 4.755282582  )}{.1cm}
\pgfcircle[fillstroke]{\pgfxy(-2.938926264 ,+ 4.045084970  )}{.1cm}
\pgfcircle[fillstroke]{\pgfxy(-4.045084972 ,+ 2.938926261 )}{.1cm}
\pgfcircle[fillstroke]{\pgfxy(-4.755282582 ,+ 1.545084968  )}{.1cm}
\pgfcircle[fillstroke]{\pgfxy(-4.755282582,- 1.545084972 )}{.1cm}
\pgfcircle[fillstroke]{\pgfxy( -4.045084972, - 2.938926262 )}{.1cm}
\pgfcircle[fillstroke]{\pgfxy(-2.938926261 ,- 4.045084972 )}{.1cm}
\pgfcircle[fillstroke]{\pgfxy(-1.545084969 ,- 4.755282582  )}{.1cm}
\pgfcircle[fillstroke]{\pgfxy(5. ,0 )}{.1cm}
\pgfcircle[fillstroke]{\pgfxy(1.545084971 ,- 4.755282582  )}{.1cm}
\pgfcircle[fillstroke]{\pgfxy(2.938926264 ,- 4.045084970  )}{.1cm}
\pgfcircle[fillstroke]{\pgfxy(4.045084972 ,- 2.938926261 )}{.1cm}
\pgfcircle[fillstroke]{\pgfxy(4.755282582 ,- 1.545084968  )}{.1cm}
\pgfellipse{\pgfxy(0,0)}{\pgfxy(6,0)}{\pgfxy(0,6)}
%Ddroite5  (Utilisateur)
\pgfsetstrokecolor{black}
\pgfsetdash{}{0pt}
\pgfxyline(0,0)(5.706339098 , 1.854101966)
\pgfxyline(0,0)( 4.854101966 ,3.526711514 )
\pgfxyline(0,0)( 3.526711513 , 4.854101966 )
\pgfxyline(0,0)(1.854101963 ,+ 5.706339099 )
\pgfxyline(0,0)(0,6)
\pgfxyline(0,0)(-1.854101965, + 5.706339098 )
\pgfxyline(0,0)(-3.526711516, + 4.854101965 )
\pgfxyline(0,0)(-4.854101966, + 3.526711513 )
\pgfxyline(0,0)(-5.706339099, + 1.854101962 )
\pgfxyline(0,0)(-6,0)
\pgfxyline(0,0)(-5.706339098 , -1.854101966)
\pgfxyline(0,0)( -4.854101966 ,-3.526711514 )
\pgfxyline(0,0)( -3.526711513 , -4.854101966 )
\pgfxyline(0,0)(-1.854101963 ,- 5.706339099 )
\pgfxyline(0,0)(0,-6)
\pgfxyline(0,0)(1.854101965, - 5.706339098 )
\pgfxyline(0,0)(3.526711516, - 4.854101965 )
\pgfxyline(0,0)(4.854101966, - 3.526711513 )
\pgfxyline(0,0)(5.706339099, - 1.854101962 )
\pgfxyline(0,0)(6,0)
               \pgfxyline(0,0)(6, 0.6306254112)
                \pgfxyline(0,0)(6, 1.275339370)
                \pgfxyline(0,0)(6, 1.949518178)
                \pgfxyline(0,0)(6, 3.464101616)
               \pgfxyline(0,0)(6, 1.275339370)
               \pgfxyline(0,0)(6, 2.671372111)
               \pgfxyline(0,0)(6, 4.359255169)
               \pgfxyline(0,0)(6,5.402424265)
                \pgfxyline(0,0)( 0.6306254112,6)
                \pgfxyline(0,0)( 1.275339370,6)
                \pgfxyline(0,0)( 1.949518178,6)
                \pgfxyline(0,0)( 3.464101616,6)
               \pgfxyline(0,0)( 1.275339370,6)
               \pgfxyline(0,0)( 2.671372111,6)
               \pgfxyline(0,0)( 4.359255169,6)
               \pgfxyline(0,0)(5.402424265,6)
                       \pgfxyline(0,0)(6, -0.6306254112)
                \pgfxyline(0,0)(6,- 1.275339370)
                \pgfxyline(0,0)(6,- 1.949518178)
                \pgfxyline(0,0)(6,- 3.464101616)
               \pgfxyline(0,0)(6, -1.275339370)
               \pgfxyline(0,0)(6, -2.671372111)
               \pgfxyline(0,0)(6, -4.359255169)
               \pgfxyline(0,0)(6,-5.402424265)
                \pgfxyline(0,0)( 0.6306254112,-6)
                \pgfxyline(0,0)( 1.275339370,-6)
                \pgfxyline(0,0)( 1.949518178,-6)
                \pgfxyline(0,0)( 3.464101616,-6)
               \pgfxyline(0,0)( 1.275339370,-6)
               \pgfxyline(0,0)( 2.671372111,-6)
               \pgfxyline(0,0)( 4.359255169,-6)
               \pgfxyline(0,0)(5.402424265,-6)
                            \pgfxyline(0,0)( -0.6306254112,-6)
                \pgfxyline(0,0)( -1.275339370,-6)
                \pgfxyline(0,0)( -1.949518178,-6)
                \pgfxyline(0,0)(- 3.464101616,-6)
               \pgfxyline(0,0)( -1.275339370,-6)
               \pgfxyline(0,0)( -2.671372111,-6)
               \pgfxyline(0,0)( -4.359255169,-6)
               \pgfxyline(0,0)(-5.402424265,-6)
                              \pgfxyline(0,0)(-6, -0.6306254112)
                \pgfxyline(0,0)(-6, -1.275339370)
                \pgfxyline(0,0)(-6, -1.949518178)
                \pgfxyline(0,0)(-6, -3.464101616)
               \pgfxyline(0,0)(-6, -1.275339370)
               \pgfxyline(0,0)(-6, -2.671372111)
               \pgfxyline(0,0)(-6, -4.359255169)
               \pgfxyline(0,0)(-6,-5.402424265)
                                           \pgfxyline(0,0)( -0.6306254112,6)
                \pgfxyline(0,0)( -1.275339370,6)
                \pgfxyline(0,0)( -1.949518178,6)
                \pgfxyline(0,0)(- 3.464101616,6)
               \pgfxyline(0,0)( -1.275339370,6)
               \pgfxyline(0,0)( -2.671372111,6)
               \pgfxyline(0,0)( -4.359255169,6)
               \pgfxyline(0,0)(-5.402424265,6)
                              \pgfxyline(0,0)(-6, 0.6306254112)
                \pgfxyline(0,0)(-6, 1.275339370)
                \pgfxyline(0,0)(-6, 1.949518178)
                \pgfxyline(0,0)(-6, 3.464101616)
               \pgfxyline(0,0)(-6, 1.275339370)
               \pgfxyline(0,0)(-6, 2.671372111)
               \pgfxyline(0,0)(-6, 4.359255169)
               \pgfxyline(0,0)(-6,5.402424265)
                              \pgfxyline(0,0)(6, 0.6306254112)
                \pgfxyline(0,0)(6, 1.275339370)
                \pgfxyline(0,0)(6, 1.949518178)
                \pgfxyline(0,0)(6, 3.464101616)
               \pgfxyline(0,0)(6, 1.275339370)
               \pgfxyline(0,0)(6, 2.671372111)
               \pgfxyline(0,0)(6, 4.359255169)
               \pgfxyline(0,0)(6,5.402424265)
                              \pgfxyline(0,0)(6, 0.6306254112)
                \pgfxyline(0,0)(6, 1.275339370)
                \pgfxyline(0,0)(6, 1.949518178)
                \pgfxyline(0,0)(6, 3.464101616)
               \pgfxyline(0,0)(6, 1.275339370)
               \pgfxyline(0,0)(6, 2.671372111)
               \pgfxyline(0,0)(6, 4.359255169)
               \pgfxyline(0,0)(6,5.402424265)
\end{tikzpicture}%
\caption{The model for the Calabi--Yau three-fold $\{x_1^{20}+x_2^6+x_3^5+x_4^4+x_5^3=0\}$ contained in $\PP(3,10,12,15,20)$.
\label{fig:bigCY}}
\end{figure}
\end{exa}

\bigskip

{\small{
\noindent \textsc{Institut Fourier, UMR du CNRS 5582,
Universit\'e de Grenoble 1,
BP 74, 38402,
Saint Martin d'H\`eres,
France}\\
\textit{E-mail address:} \url{chiodo@ujf-grenoble.fr}

\vspace{.3cm}

\noindent \textsc{Department of Mathematics, University of Michigan, Ann Arbor, MI 48109-1109,
USA} and  \textsc{Yangtze Center of Mathematics, Sichuan University, Chengdu,
610064, P.R. China}\\
\textit{E-mail address:} \url{ruan@umich.edu}
}
}
\end{document}